\documentclass[MSNbibl,number,citesort,dvips]{arxbj}
\usepackage{mathbh,mathrsfs}
\usepackage{graphicx}
\usepackage{url, breakurl}
\aid{0}
\volume{20}
\issue{3}
\pubyear{2014}
\firstpage{1260}
\lastpage{1291}
\doi{10.3150/13-BEJ521} 

\makeatletter
\newcommand{\bigtimes}{\mathop{\!\mbox{\parbox[c][9pt][b]{18pt}{\fontsize{18}{18}\selectfont{$\times$}}}\!\!\!\!}}

\def\bbm{\bolds}
\def\bm{\mathbf}
\newcommand{\eqref}[1]{(\ref{#1})}
\newcommand{\indic}{\mathbh{1}}
\newcommand{\X}{\mathbb{X}}
\newcommand{\Y}{\mathbb{Y}}

\newcommand{\R}{\mathbb{R}}
\newcommand{\E}{\mathbb{E}}
\renewcommand{\P}{\mathbb{P}}

\newcommand{\Fcr}{\mathscr{F}}
\newcommand{\Dcr}{\mathscr{D}}
\newcommand{\Xcr}{\mathscr{X}}
\newcommand{\Mcr}{\mathscr{M}}

\newcommand{\Xf}{\mathfrak{X}}

\newcommand{\ddr}{\mathrm{d}}
\newcommand{\edr}{\mathrm{e}}

\newcommand{\Lc}{\mathcal{L}}

\newcommand{\ep}{\varepsilon}
\newcommand{\wt}{\widetilde}

\def\simind{\mathop{\sim}\limits^{\mathrm{ind}}}
\def\simiid{\mathop{\sim}\limits^{\mathrm{i.i.d.}}}

\newcommand{\mt}{\tilde{\mu}} 
\newcommand{\pt}{\tilde{p}} 
\newcommand{\La}{\mathcal{L}}
\newcommand{\Ac}{\mathcal{A}}
\newtheorem{proposition}{Proposition}
\newtheorem{coro}{Corollary}
\newproclaim{definition}{Definition}
\newremark{example}{Example}
\newremark{remark}{Remark}
\makeatother

\begin{document}
\begin{frontmatter}

\title{Bayesian inference with dependent normalized completely random measures}
\runtitle{Dependent normalized CRMs}

\begin{aug}
\author[a,c]{\inits{A.}\fnms{Antonio}~\snm{Lijoi}\corref{}\thanksref{a,c,e1}\ead[label=e1,mark]{lijoi@unipv.it}},
\author[b,c]{\inits{B.}\fnms{Bernardo} \snm{Nipoti}\thanksref{b,c,e2}\ead[label=e2,mark]{bernardo.nipoti@unito.it}}
\and
\author[b,c]{\inits{I.}\fnms{Igor} \snm{Pr\"unster}\thanksref{b,c,e3}\ead[label=e3,mark]{igor.pruenster@unito.it}}
\runauthor{A. Lijoi, B. Nipoti and I. Pr\"unster} 
\address[a]{Department of Economics and Management, University of
Pavia, Via San
Felice 5, 27100 Pavia, Italy. \printead{e1}}
\address[b]{Department of Economics and Statistics, University of
Torino, Corso Unione Sovietica 218/bis, 10134 Torino, Italy. \printead{e2,e3}}
\address[c]{Collegio Carlo Alberto, Via Real Collegio 30, 10024
Moncalieri, Italy}
\end{aug}

\received{\smonth{3} \syear{2012}}
\revised{\smonth{10} \syear{2012}}
%
\begin{abstract}
The proposal and study of dependent prior processes
has been a
major research focus in the recent Bayesian nonparametric literature.
In this paper, we introduce a flexible class of dependent nonparametric
priors, investigate their properties and derive a suitable sampling
scheme which allows their concrete implementation. The proposed class
is obtained by normalizing dependent completely random measures, where
the dependence arises by virtue of a suitable construction of the
Poisson random measures underlying the completely random measures. We
first provide general distributional results for the whole class of
dependent completely random measures and then we specialize them to two
specific priors, which represent the natural candidates for concrete
implementation due to their analytic tractability: the bivariate
Dirichlet and normalized $\sigma$-stable processes. Our analytical
results, and in particular the partially exchangeable partition
probability function, form also the basis for the determination of a
Markov Chain Monte Carlo algorithm for drawing posterior inferences,
which reduces to the well-known Blackwell--MacQueen P\'olya urn scheme
in the univariate case. Such an algorithm can be used for density
estimation and for analyzing the clustering structure of the data and
is illustrated through a real two-sample dataset example.
\end{abstract}

%
\begin{keyword}
\kwd{completely random measure}
\kwd{dependent Poisson processes}
\kwd{Dirichlet process}
\kwd{generalized Pol\'ya urn scheme}
\kwd{infinitely divisible vector}
\kwd{normalized $\sigma$-stable process}
\kwd{partially exchangeable random partition}
\end{keyword}

\end{frontmatter}

\section{Introduction}\label{sec1}

The construction of dependent random probability measures for Bayesian
inference has attracted considerable attention in the last decade. The
seminal contributions of MacEachern \cite{MacEachern,MacEachernTR},
who introduced a general class of dependent processes including a
popular dependent version of the Dirichlet process, paved the way to a
burst in the literature on (covariate) dependent processes and their
application in a variety of frameworks such as, for example,
nonparametric regression, inference on time series data, meta-analysis,
two-sample problems. Reviews and key references can be found in, for
example, \cite{mq,dunsone,tehjord}. Most contributions to this line of
research rely on random probability measures defined by means of a
stick-breaking procedure, a popular method set forth in its generality
for the first time in \cite{ilancillotto}. Dependence among different
\textit{stick-breaking priors} is created by indexing either the
stick-breaking weights or the locations or both to relevant
covariates. To be more specific, if $\mathcal{Z}$ denotes the covariate
space and $\{(\omega_{j,z})_{j\ge1}\dvt z\in\mathcal{Z}\}$ is a
collection of sequences of independent nonnegative weights, the
stick-breaking procedure consists in defining $ p_{1,z}=\omega_{1,z}$
and $p_{j,z}=\omega_{j,z} \prod_{i=1}^{j-1} (1-\omega_{i,z})$. A
typical choice is\vspace*{1pt} then $\omega_{i,z}\sim\operatorname{Beta}(a_{i,z},b_{i,z})$
with parameters $(a_{i,z},b_{i,z})$ such that\vspace*{1pt} $\sum_{j\ge1}p_{j,z}=1$,
almost surely. If one further considers collections of sequences $\{
(X_{i,z})_{i\ge1}\dvt z\in\mathcal{Z}\}$ with the $X_{i,z}$, for $i\ge
1$, taking values in a space $\X$ and i.i.d. from a nonatomic
probability measure $P_{0,z}$, a \textit{covariate dependent random
probability measure} $\tilde p_z=\sum_{j\ge1} p_{j,z} \delta
_{X_{j,z}}$ is obtained. The dependence between weights $\omega_{i,z}$
and $\omega_{j,z'}$ and/or between the support points $X_{i,z}$ and
$X_{j,z'}$, for $z\ne z'$, induces dependence between $\tilde p_z$ and~$\tilde p_{z'}$. This general framework is then tailored to the
specific application at issue. One of the main reasons of the success
of stick-breaking constructions is their attractiveness from
computational point of view along with their flexibility since, as
shown in \cite{bjq}, they have full weak support under mild
assumptions. On the other hand, a drawback is represented by the
difficulty of studying their distributional properties due to their
analytical intractability.
In this paper, we propose a radically different approach to the
construction of dependent nonparametric priors that relies on \textit
{completely random measures} (CRMs) introduced by Kingman \cite
{kingm93}. For the case of exchangeable setting, in \cite{lp10} it has
been shown that CRMs represent a unifying concept of the Bayesian
Nonparametrics given most discrete nonparametric priors can be seen as
transformations of CRMs. Our general plan consists in defining a broad
class of dependent CRMs thus obtaining a vector of dependent random
probability measures via a suitable transformation. A relevant
motivation for undertaking such an approach is represented by the
consideration that the study of distributional properties of the models
are essential for their deep understanding and sound applications. In
this respect, even though CRMs are infinite-dimensional objects, they
can be summarized by a single measure, that is, their intensity, which
allows to derive key distributional properties.

\subsection{Dependent Poisson random measures}\label{sec1.1} A key idea of our
approach consists in defining dependent CRMs by creating dependence at
the level of the underlying Poisson random measures (PRM). To this end,
we resort to a class of bivariate dependent PRMs devised by Griffiths
and Milne in \cite{GM78}. In particular, let $\tilde N$ be a PRM on
$\Y
$ with intensity measure $\bar\nu$. The corresponding Laplace
functional transform, which completely characterizes the PRM, is then
given by
\[
\E \bigl[\edr^{-\int f\, \ddr\tilde N} \bigr]=\edr^{-\bar\nu
(1-\edr^{-f})}
\]
for any measurable function $f\dvtx \Y\to\R$ such that $\int|f| \,\ddr
\tilde N<\infty$ (a.s.). Recall also that a Cox process is a PRM with
random intensity. See \cite{daley} for an exhaustive account. Consider
now a vector of (possibly dependent) PRMs $(\tilde N_1,\tilde N_2)$ on
$\Y$ with the same marginal intensity measure~$\bar\nu$. Griffiths and
Milne \cite{GM78} prove that the $\tilde N_i$'s admit an additive
representation
%
%
\begin{equation}
\label{eqprm} \tilde N_i= M_i+M_0,\qquad i=1,2,
\end{equation}
where $M_1$, $M_2$ and $M_0$ are independent Cox processes with
respective random intensities $\nu$, $\nu$ and $\nu_0$ such that
$\nu
_0\le\bar\nu$ (a.s.) and $\nu=\bar\nu-\nu_0$ if and only if the Laplace
transform has the following form
%
%
\begin{equation}
\label{eqprmjointlapl} \E \bigl[\edr^{-\tilde N_1(f_1)-\tilde N_2(f_2)} \bigr] =\edr^{-\sum_{i=1}^2\bar\nu(1-\edr^{-f_i})}
\varphi\bigl(\bigl(1-\edr^{-f_1}\bigr) \bigl(1-\edr^{-f_2}\bigr)
\bigr)\vadjust{\goodbreak}
\end{equation}
for some functional $\varphi$. Such a result is appealing for at least
two reasons. From an intuition point of view, it provides a neat
additive representation \eqref{eqprm} of the $\tilde N_i$'s with a
common and idiosyncratic component, $M_0$ and $M_i$, for $i=1,2$,
respectively. From an operational point of view, it yields a well
identified structure \eqref{eqprmjointlapl} for the Laplace
functional, which becomes completely explicit in the cases where one is
able to determine the form of $\varphi$. In fact, when working with PRMs
and CRMs, the Laplace functional is the main operational tool for
deriving analytical results useful for Bayesian inference and such a
relatively simple structure is actually quite surprising for the
dependent case.

The pair of PRMs constructed according to \eqref{eqprm} is, then, used
to define a vector of dependent CRMs $(\tilde\mu_1,\tilde\mu_2)$.
Recall that CRMs are random measures giving rise to mutually
independent random variables when evaluated on pairwise disjoint
measurable sets. Moreover, they can always be represented as
functionals of an underlying PRM, which in the particular case of $\Y
=\R
^+ \times\R^d$ corresponds to the celebrated \textit{L\'evy--Ito
decomposition}. Therefore, by setting $\Y=\R^+\times\X$, from
$(\tilde
N_1,\tilde N_2)$ one can define the corresponding vector of CRMs $(\mt
_1,\mt_2)$ with components given by $\mt_i(\ddr x)=\int_{\R^+}s
\tilde
N_i(\ddr s,\ddr x)$.

Finally, a vector of dependent random probability measures on $\X$ is
obtained as $(\tilde p_1,\tilde p_2)\stackrel{d}{=}(T(\mt_1),T(\mt_2))$
where $T$ is a transformation of the CRM such that $T(\mt_i)(\X)=1$
a.s. Here we focus on one of the most intuitive transformations, namely
``normalization'', which corresponds to $T(\mt)=\mt/\mt(\X)$. Such a
normalization procedure is widely used in the univariate case. Already
Ferguson \cite{ferg} showed that the Dirichlet process can be defined
as normalization of a gamma CRM. Such a procedure has then been
extended and analyzed for general univariate CRMs in \cite
{rlp,jlp06,jlp}. More recently, an interesting construction of a
subclass of normalized CRMs has been proposed in \cite{orbanzo}. See
\cite{lp10} for
a review of other commonly used transformations~$T$.

In the literature there are already some proposals, although not in a
general framework and analytical depth as set forth here, making use of
dependent CRMs for defining dependent random probability measures. For
example, in \cite{cipriota} and in \cite{raoteh} one can find a model
that coincides with a special case we consider in this paper, namely a
version of the bivariate Dirichlet process. In these two papers, the
authors devise samplers that take advantage of a mixture representation
of $\tilde p_1$ and of $\tilde p_2$ whose weights are, only for their
special case, independent from the $\tilde p_i$'s. In a similar
fashion, \cite{mqr04} proposes dependent convex linear combinations of
Dirichlet processes as a tool for examining data originated from
different experiments. Vector CRMs, whose dependence is induced by
suitable L\'evy copulas, are proposed in \cite{ilenia} for defining a
vector of dependent neutral to the right processes and in \cite
{fabrizio} in order to introduce a bivariate two-parameter
Poisson--Dirichlet process. In addition to the great generality of our
results, two important features of our proposal are to be highlighted:
it preserves computational efficiency since we are able to deduce a
generalization of the Blackwell--MacQueen urn scheme for the dependent
setting implementable in real-world applications, and it sheds light
on theoretical properties of the vector of random probability measures
we are proposing, therefore improving the understanding of the model.

\subsection{Goals and outline of the paper}\label{sec1.2}

As mentioned above, we will investigate vectors of random probabilities
$(\pt_1,\pt_2)$ obtained by normalizing pairs of dependent CRMs
$(\tilde
\mu_1,\tilde\mu_2)$. The distribution of $(\pt_1,\pt_2)$ plays the role
of mixing measure in the representation of the law of a pair of
partially exchangeable sequences or, in other terms, of prior
distribution for a partially-exchangeable observation process. We will
determine an expression for the probability distribution of the
partially exchangeable partition induced by $(\tilde p_1,\tilde p_2)$.
Such a result will also lead us to achieve an extension of the
univariate Blackwell--MacQueen P\'olya urn scheme. The corresponding
Gibbs sampler is then implemented to draw a full Bayesian analysis for
density estimation and cluster analysis in two-sample problems. The
general results will, then, be specialized to two specific priors
where: (i) the $\tilde\mu_i$'s are gamma CRMs thus yielding a vector of
dependent Dirichlet processes; (ii) the $\tilde\mu_i$'s are $\sigma
$-stable CRMs that give rise to a vector of dependent normalized
$\sigma$-stable processes.

The outline of the paper is as follows. In Section~\ref{sec2}, we introduce some
notation and formalize the form of dependence we briefly touched upon
before. In Section~\ref{secpeppf}, we consider pairs of partially exchangeable
sequences directed by the distribution of $(\pt_1,\pt_2)$ and describe
some of their distributional properties. Section~\ref{secddp} considers dependent
mixtures and introduces the main distributional tools that are needed
for their application to the analysis of partially exchangeable data.
Section~\ref{sec5} provides a description of the prior specification we adopt
and the sampler we resort to. Finally, Section~\ref{sec6} contains an
illustration with a real dataset which is analyzed through mixture
models with both dependent Dirichlet and normalized $\sigma$-stable.
The proofs are postponed to the \hyperref[app]{Appendix}. A~key tool for proving our
results is represented by an extension to the partial exchangeable case
of a technique introduced and subsequently refined in \cite
{ip,jlp06,jlp}. Such a technique was originally developed for deriving
conditional distributions of normalized random measures \cite{rlp} but,
as highlighted in \cite{lp10}, it can be actually applied to any
exchangeable model based on completely random measures. Therefore, it
is worth remarking that the extension to the partial exchangeable setup
is also of independent interest.\looseness=1

\section{Dependent completely random measures}\label{sec2}
Let us start by stating more precisely some of the concepts sketched in
the \hyperref[sec1]{Introduction}. Consider a probability space $(\Omega,\Fcr,\P)$ and
denote by $M_\X$ the set of boundedly finite measures on a complete and
separable metric space $\X$. Further, the Borel $\sigma$-algebras on
$M_\X$ and $\X$ are denoted by $\Mcr_\X$ and $\Xf$, respectively. A
completely random measure (CRM) $\mu$ on $(\X,\Xf)$ is a measurable
function on $(\Omega, \Fcr,\P)$ taking values in $(M_\X,\Mcr_\X)$
such that for any $ d \geq  1$ and any collection
$\{ A_1,\ldots,A_d\}$ of pairwise disjoint sets in $\mathfrak{X}$, the random variables
$\mu (A_1),\ldots,\mu(A_d)$ are mutually independent.
It is well known that
if $N$ is a Poisson random measure on $\R^+\times\X$, then
%
%
\begin{equation}
\label{eqpoissfunct} \mu(B)=\int_{\R^+\times B} s N(\ddr s,\ddr x)\qquad
\forall B\in\Xf
\end{equation}
is a CRM on $(\X,\Xf)$. See \cite{kingm93,daley} and, for example, \cite
{james05} for uses of representation \eqref{eqpoissfunct} for
Bayesian modeling. If $\bar\nu$ is the intensity of $N$ and for brevity
$\mu(f):=\int f \,\ddr\mu$, the Laplace exponent of $\mu(f)$ is of
the form
%
%
\begin{equation}
\label{eqlapl} -\log \bigl(\E \bigl[\edr^{-\mu(f)} \bigr] \bigr)= \int
_{\R^+\times\X}\bigl[1-\edr^{-sf(x)}\bigr] \bar\nu(\ddr s,\ddr
x)=:\psi(f)
\end{equation}
for any measurable function $f\dvtx \X\to\R$ such that $\mu(|f|)=\int|f|
\,\ddr\mu<\infty$, almost surely. By virtue of \eqref
{eqpoissfunct}, we
can construct dependent CRMs as linear functionals of dependent PRMs
determined according to \eqref{eqprm}. To state it more precisely, let
$P_0$ be a nonatomic probability measure on $(\X,\Xf)$ and $r(\ddr
s)=\rho(s) \,\ddr s$ a (possibly infinite) measure on $\R^+$. Suppose,
further, that $\tilde N_1$ and $\tilde N_2$ are defined as in \eqref
{eqprm}, where $M_1$, $M_2$ and $M_0$ are three independent Cox
processes with respective random intensities $\nu$, $\nu$ and $\nu_0$
such that $\nu+\nu_0=\bar\nu$, almost surely. Henceforth, we shall
assume $\bar{\nu}(\ddr s,\ddr x)=cP_0(\ddr x) \rho(s) \,\ddr s$.

\begin{definition}\label{def1}
Let $(\tilde N_1, \tilde N_2)$ be a
vector of Griffiths--Milne (GM) dependent PRMs as in \eqref{eqprm} and
define the CRMs $\tilde\mu_i(\ddr x)=\int_{\R^+}s \tilde N_i(\ddr
s,\ddr x)$, for $i=1,2$. Then $(\tilde\mu_1,\tilde\mu_2)$ is said
to be
a vector of \textit{GM-dependent} CRMs. The marginal intensity of
$\tilde\mu_i$ coincides with $\bar\nu$.
\end{definition}

In the sequel, we will focus on a simple class of Cox
processes defined through an intensity of the form
%
%
\begin{equation}
\label{eqintensity} \nu(\ddr s,\ddr x) = cZ P_0(\ddr x) \rho(s) \,
\ddr s %
\end{equation}
for some $[0,1]$-valued random variable $Z$.
To ease the exposition, and with no loss of generality, we will work
conditionally on a fixed value $Z=z$ which makes the Cox processes in
\eqref{eqprm} coincide with PRMs. According to the definition above,
the marginals of a vector of GM-dependent CRMs are equally distributed and
%
%
\begin{equation}
\label{eqdecomp} \mt_i (\ddr x) = \int_{\R^+}s
M_i(\ddr s,\ddr x)+\int_{\R^+}s M_0(
\ddr s,\ddr x) =\mu_i(\ddr x)+\mu_0(\ddr x),
\end{equation}
where $\mu_i$, with $i=1,2$, and $\mu_0$ are independent CRMs with
Laplace functional transforms
\[
\E \bigl[\edr^{-\mu_i(f)} \bigr]=\edr^{-cz\psi(f)},\qquad \E \bigl[
\edr^{-\mu_0(f)} \bigr]=\edr^{-c(1-z)\psi(f)},
\]
%
where $\psi$ is defined as in \eqref{eqlapl}. Given the simple form of
the intensities specified in \eqref{eqintensity}, one can determine
the form of $\varphi$ in \eqref{eqprmjointlapl} explicitly and
straightforwardly obtains a tractable expression for the joint Laplace
functional transform of $(\mt_1,\mt_2)$ given by
%
%
\begin{equation}
\label{eqjointlapl} \E \bigl[\edr^{-\mt_1(f_1)-\mt_2(f_2)} \bigr]=\edr^{-cz[\psi
(f_1)+\psi
(f_2)]-c(1-z)\psi(f_1+f_2)}
\end{equation}
for any pair of measurable functions $f_i\dvtx \X\to\R$, for $i=1,2$, such
that $\P[\mt_i(|f_i|)<\infty]=1$. In order to further clarify the above
concepts and construction, let us consider two special cases involving
well-known CRMs.

\begin{example}[(Gamma process)]\label{ex1}
Set $\rho
(s)=\edr
^{-s} s^{-1}$ in \eqref{eqintensity} which results in $\mu$ being a
gamma CRM.
The corresponding Laplace exponent reduces to $\psi (f) =\int \log(1 +
f )\,\mathrm{d} P_0$ for any
measurable function $f$ such that $\int\log(1 + |f|)\,\mathrm{d}P_0 <\infty$.
If $f_i\dvtx  \X\to\R$
are, for
$i=1,2$, measurable functions such
that $\int \log(1 +|f_i|)\,\mathrm{d}P_0 <\infty$, one has
%
\[
\E \bigl[\edr^{-\mt_1(f_1)-\mt_2(f_2)} \bigr]=\edr^{-c\int\log
(1+f_1+f_2)\,\ddr P_0-
cz\int\log{(1+f_1)(1+f_2)}/{(1+f_1+f_2)} \,\ddr P_0}.
\]
\end{example}

\begin{example}[($\boldsymbol{\sigma}$-stable process)]\label{ex2}
Set
$\rho(s)=\sigma s^{-1-\sigma}/\Gamma(1-\sigma)$, with $\sigma\in
(0,1)$, in \eqref{eqintensity}
which results in $\mu$ being a $\sigma$-stable CRM.
The corresponding
Laplace exponent
reduces to $\psi (f)= \int f^{\sigma} \,\mathrm{d}P_0$ for any
measurable function
$f$ such that $\int|f|^{\sigma} \,\mathrm{d}P_0 <\infty$.
Let $f_i\dvtx  \X\to\R$
be such that $\int|f_i|^\sigma
\,\ddr
P_0<\infty$, for $i=1,2$. Then
\[
\E \bigl[\edr^{-\mt_1(f_1)-\mt_2(f_2)} \bigr]=\edr^{-cz\int
(f_1^\sigma
+f_2^\sigma)\,\ddr P_0-
c(1-z)\int(f_1+f_2)^\sigma\,\ddr P_0}.
\]
\end{example}

The final step needed for obtaining the desired vector of dependent
random probability measures consists in normalizing the previously
constructed CRMs, in the same spirit as in \cite{rlp} for the
univariate case. To perform the normalization, we need to ensure $\P
[\mu
_i(\X)\in(0,\infty)]=1$, for $i=0,1,2$, which is guaranteed by
requesting $\int_0^\infty\rho(s) \,\ddr s=\infty$ (see \cite{rlp}) and
corresponds to considering CRMs which jump infinitely often on any
bounded set. By normalizing $\mt_1$ and $\mt_2$, we can then define the
vector of dependent random probability measures
%
%
\begin{equation}
(\tilde p_1, \tilde p_2 )\stackrel{d} {=} \bigl(
\mt_1/\mt_1(\X), \mt_2/\mt_2(\X)
\bigr) \label{eqvectornrmi}
\end{equation}
to be termed \textit{GM-dependent normalized CRM} in the following.

Having described the main concepts and tools we are resorting to, our
next goal is the application of $(\pt_1,\pt_2)$ as a nonparametric
prior for the statistical analysis of partially exchangeable data.

\section{Partially exchangeable sequences}\label{secpeppf}
For our purposes, we resort to the notion of partial exchangeability as
set forth by de~Finetti in~\cite{defin38} and described as follows. Let
$X=(X_n)_{n\ge1}$ and $Y=(Y_n)_{n\ge1}$ be two sequences of $\X
$-valued random elements defined on some probability space $(\Omega
,\Fcr,\P)$ and $P_\X$ is the space of probability measures on $(\X
,\Xf
)$. If $\bm{X}^{(n_1)}=(X_1,\ldots,X_{n_1})$ and $\bm
{Y}^{(n_2)}=(Y_1,\ldots,Y_{n_2})$ are the first $n_1$ and $n_2$ values
of the sequences $X$ and $Y$, respectively, we have
%
%
\begin{equation}
\label{eqpartialexchange} \P \bigl[\bm{X}^{(n_1)}\in A_1,
\bm{Y}^{(n_2)}\in A_2 \bigr] =\int_{P_\X^2}p_1^{n_1}(A_1)
p_2^{n_2}(A_2) \Phi(\ddr p_1,\ddr
p_2)
\end{equation}
for any $A_1\in\Xf^{n_1}$, $A_2\in\Xf^{n_2}$, with $p_i^{n_i}$ being
the $n$-fold product measure $p_i\times\cdots\times p_i$ and $\Phi
$ is a probability distribution on $P_\X^2=P_\X\times P_\X$ which acts
as nonparametric prior for Bayesian inference. We also denote as $\Phi
_i$ the marginal distribution of $\tilde p_i$ on $P_\X$. Since $\tilde
p_i$ is a normalized CRM, then the weak support of $\Phi_i$ contains
all probability measures on $\X$ whose support is contained in the
support of the base measure $P_0$. Hence, if the support of $P_0$
coincides with $\X$, a GM-dependent normalized CRM $(\tilde
p_1,\tilde
p_2)$ has full weak support with respect to the product topology on
$P_\X^2$. Having a large support is a minimal requirement a
nonparametric prior must comply with in order to ensure some degree of
flexibility in statistical analysis.

It should be also noted that the dependence structure displayed in
assumption~\eqref{eqpartialexchange} is also the starting point in
\cite{cifreg78} where the authors propose an example (the first we are
aware of in the literature) of nonparametric prior for partially
exchangeable arrays which coincides with a mixture of products of
Dirichlet processes. Furthermore, \eqref{eqpartialexchange} defines
the framework in which recent proposals of dependent nonparametric
priors can be embedded.

\subsection{Dependence between \texorpdfstring{$\tilde p_1$ and $\tilde p_2$}{p1 and p2}}\label{sec3.1}
An important preliminary result we state concerns the mixed moment of
$(\tilde p_1(A),\tilde p_2(B))$ for any $A$ and $B$ in $\Xf$. To this
end, define the following quantity\vspace*{-2pt}
%
%
\begin{equation}
\label{eqtau} \tau_q(u):=\int_0^\infty
s^q \edr^{-us} \rho(s) \,\ddr s
\end{equation}
for any $q\ge1$. Moreover, to simplify the notation in (\ref{eqlapl}) we set $\psi(u)=\psi
(u\indic_\X)$ for any $u>0$, where $\indic_A$ is the indicator function
on set $A$. One can, then, prove the following proposition.\vspace*{-2pt}

\begin{proposition}\label{pr1}
Let $(\tilde p_1, \tilde p_2)$ be a
vector of GM-dependent normalized CRM defined in \eqref
{eqvectornrmi}. For any $A$ and $B$ in $\Xf$ one has
%
%
\begin{eqnarray}
\label{eqmixmom} \E \bigl[\tilde p_1(A) \tilde p_2(B)
\bigr]&=&P_0(A)P_0(B)+ \bigl[P_0(A\cap B)
-P_0(A)P_0(B) \bigr]
\nonumber
\\[-8pt]
\\[-8pt]
\nonumber
&&{}\times c(1-z) \int_{(0,\infty)^2}\edr^{-cz[\psi(u)+\psi(v)]-c(1-z)\psi
(u+v)}\tau
_2(u+v) \,\ddr u \,\ddr v.
\end{eqnarray}
Moreover, it follows that\vspace*{-2pt}
%
%
\begin{equation}
\label{eqcorrelation} \operatorname{Corr}\bigl(\tilde p_1(A),\tilde
p_2(B)\bigr)= \frac{(1-z) [P_0(A\cap B)-P_0(A)P_0(B) ]}{\sqrt{P_0(A)[1-P_0(A)]}
\sqrt{P_0(B)[1-P_0(B)]}} \mathcal{I}(c,z),
\end{equation}
where\vspace*{-2pt}
\[
\mathcal{I}(c,z):= \frac{\int_0^\infty\int_0^\infty\edr
^{-cz[\psi
(u)+\psi(v)]-c(1-z)\psi(u+v)}
\tau_2(u+v) \,\ddr u\,\ddr v}{\int_0^\infty u \edr^{-c\psi(u)}
\tau_2(u) \,\ddr u}.
\]
\end{proposition}

It can be easily seen that if $A=B$, then the correlation in
\eqref{eqcorrelation} reduces to $(1-z)\mathcal{I}(c,z)$ and does not
depend on the specific set where the two random probabilities $\tilde
p_1$ and $\tilde p_2$ are evaluated. This fact is typically used to
motivate $(1-z)\mathcal{I}(c,z)$ as a measure of the (overall)
dependence between $\tilde p_1$ and $\tilde p_2$. Coherently with our
construction $\tilde p_1$ and $\tilde p_2$ are uncorrelated if $z=1$,
and the same can be said if $A$ and $B$ are independent with respect to
the baseline probability measure $P_0$. The previous expression is
structurally neat and, as will be shown in the following illustrations,
in some important special cases the double integral $\mathcal{I}(c,z)$
can be made sufficiently explicit so to allow a straightforward computation.\vspace*{-2pt}

\setcounter{example}{0}
\begin{example}[(Continued)]
If $\mt_1$, $\mt_2$ are two
dependent CRMs, one has $\tau_q(u)=\Gamma(q) (1+u)^{-q}$ and the
correlation between the corresponding GM-dependent Dirichlet processes
coincides with~\eqref{eqcorrelation} where\vspace*{-2pt}
%
%
\begin{equation}
\mathcal{I}(c,z)=\frac{c}{c+1} {}_3F_2(c-cz+2,1,1;c+2,c+2;1),
\label{eqgenhypergeom}
\end{equation}
where ${}_3F_2$ is the generalized hypergeometric function\vspace*{-2pt}
%
%
\begin{equation}
\label{def3F2} {}_3F_2(\alpha,\beta,\rho;\gamma,\sigma;x)=
\sum_{j\ge0}\frac
{(\alpha
)_j(\beta)_j(\rho)_j}{
j!(\gamma)_j(\sigma)_j} x^j
\end{equation}
and $(a)_n = \Gamma(a + n)/ \Gamma(a+n)$ for any\vadjust{\goodbreak} $a> 0$ and
any non-negative integer $n$.
The above series converges if $|x|<1$ and it does for $x=1$ provided
that Re$(\gamma+\sigma-\alpha-\beta-\rho)>0$, with Re$(z)$
denoting the
real part of a complex number $z$.
\end{example}

\begin{example}[(Continued)]
If $\mt_1$, $\mt_2$ are
$\sigma$-stable dependent CRMs, one has $\tau_q(u)=\sigma(1-\sigma
)_{q-1} u^{\sigma-q}$ and the correlation between the corresponding
dependent normalized $\sigma$-stable processes is equal to \eqref
{eqcorrelation} with\vspace*{-1pt}
\[
\mathcal{I}(c,z)=\frac{1}{\sigma} \int_0^1
\frac{w^{1/\sigma-1}}{[1+z(1-w^{1/\sigma})^\sigma-z(1-w)]} \,\ddr w.
\]
Even if we are not able to evaluate the above integral analytically, a
numerical approximation can be easily determined.
\end{example}

\subsection{Partition probability function}\label{secpepf}\vspace*{-3pt} The
procedure adopted for determining an expression for the mixed moments
of $\tilde p_1$ and $\tilde p_2$ can be extended to provide a form for
the partially exchangeable partition probability function (pEPPF) for
the $n_1+n_2$ random variables (r.v.'s)
$\bm{X}^{(n_1)}$ and $\bm{Y}^{(n_2)}$. It is worth recalling that the
concept of EPPF plays an important role in modern probability theory
(see \cite{pitman06} and references therein) and, implicitly, in
numerous MCMC algorithms one ends up ``sampling from the partition'' as
well. First, note that if $z<1$\vspace*{-1pt}
\[
\P[X_i=Y_j]=c\int_0^\infty
u \edr^{-c\psi(u)} \tau_2(u) \,\ddr u>0
\]
for any $i$ and $j$: hence, with positive probability any of the
elements of the first sample $\bm{X}^{(n_1)}$ can coincide with any
element from $\bm{Y}^{(n_2)}$.
This leads us to address the issue of determining the probability that
the two samples are partitioned into $K=K_{1}+K_{2}+K_{0}$ clusters of
distinct values where
\begin{longlist}[(a)]
\item[(a)] $K_1$ is the number of distinct values in the first sample
$\bm{X}^{(n_1)}$ not coinciding with any of the $Y_j$'s;
\item[(b)] $K_2$ is the number of distinct values in the second sample
$\bm{Y}^{(n_2)}$
not coinciding with any of the $X_j$'s;
\item[(c)] $K_0$ is the number of distinct values that are shared by
both samples
$\bm{X}^{(n_1)}$ and~$\bm{Y}^{(n_2)}$.
\end{longlist}
Moreover, we denote by $\bm{N}^{(i)}=(N_{1,i},\ldots,N_{K_i,i})$ the
vector of frequencies for the $K_i$ unshared clusters and with $\bm
{Q}^{(i)}=(Q_{1,i},\ldots,Q_{K_0,i})$ the vector of frequencies the
sample $\bm{X}^{(n_1)}$, if $i=1$, or the sample $\bm{Y}^{(n_2)}$, if
$i=2$, contributes to each of the shared clusters. Correspondingly, we
introduce the sets of vectors of positive integers\vspace*{-1pt}
\[
\Delta_{n_i,k_i,k_0}:= \Biggl\{\bigl(\bm{n}^{(i)},
\bm{q}^{(i)}\bigr)\dvt \sum_{l=1}^{k_i}n_{l,i}+
\sum_{r=1}^{k_0} q_{r,i}=n_i
\Biggr\},
\]
where the more concise notation $
\bm{n}^{(i)}=(n_{1,i},\ldots,n_{k_i,i})$ and $
\bm{q}^{(i)}=(q_{1,i},\ldots,q_{k_0,i})$ is used,
for $i=1,2$. The result\vadjust{\goodbreak}
we are going to state characterizes the probability distribution
of the random partition induced by $(\mathbf{X}^{(n_1)}, \mathbf{Y}^{(n_2)})$
as
encoded
by the vector of positive integers
$(K_1,K_2,K_0,\bm{N}^{(1)},\bm{N}^{(2)},\bm{Q}^{(1)},\bm{Q}^{(2)})$.
Such a distribution has masses at points $(k_1,k_2,k_0,\bm
{n}^{(1)}, \bm
{n}^{(2)},\bm{q}^{(1)},\bm{q}^{(2)})$ that we denote as
$\Pi_{k}^{(n_1+n_2)}(\bm{n}^{(1)},\bm{n}^{(2)},\bm{q}^{(1)},\bm
{q}^{(2)})$, where $k=k_1+k_2+k_0$.

\begin{proposition}\label{pr2} Let $(\tilde p_1, \tilde
p_2)$ be a GM-dependent normalized CRM defined in \eqref
{eqvectornrmi}. For any $(\bm{n}^{(i)},\bm{q}^{(i)})\in\Delta
_{n_i,k_i,k_0}$, with $i=1,2$, and for any nonnegative integers $k_1$,
$k_2$ and\vadjust{\eject} $k_0$ such that $k_l+k_0\in\{1,\ldots,n_l\}$, for $l=1,2$,
one has
\begin{eqnarray*}
&&\Pi_{k}^{(n_1+n_2)}\bigl(\bm{n}^{(1)},
\bm{n}^{(2)},\bm{q}^{(1)},\bm{q}^{(2)}\bigr)\\
&&\quad=
\frac{c^{k}}{\Gamma(n_1)\Gamma(n_2)} \sum_{(*)}(1-z)^{k_0+|\bm{i}|+|\bm{l}|}
z^{k_1+k_2-|\bm{i}|-|\bm{l}|}
\\
&&\qquad{}\times \int_0^\infty\int_0^\infty
u^{n_1-1} v^{n_2-1} \edr ^{-cz[\psi
(u)+\psi(v)]
-c(1-z)\psi(u+v)}
\\
&&\hspace*{49pt}\qquad{}\times\prod
_{j=1}^{k_1}\tau_{n_{j,1}}(u+i_j v)\prod_{j=1}^{k_2}
\tau_{n_{j,2}}(l_j u+ v) \prod_{r=1}^{k_0}
\tau_{q_{r,1}+q_{r,2}}(u+v) \,\ddr u \,\ddr v,
\end{eqnarray*}
where the sum runs over the set of all vectors of integers $\bm
{i}=(i_1,\ldots,i_{k_1})\in
\{0,1\}^{k_1}$ and $\bm{l}=(l_1,\ldots,l_{k_2})\in\{0,1\}^{k_2}$,
whereas $|\bm{i}|=\sum_{j=1}^{k_1} i_j$ and $|\bm{l}|=\sum_{j=1}^{k_2} l_j$.
\end{proposition}

The expression, though in closed form and of significant
theoretical interest, is quite difficult to evaluate due to the
presence of the sum with respect to the integer vectors $\bm{i}$ and
$\bm{l}$. Nonetheless, Proposition~\ref{pr2} is going to be a fundamental tool
for the derivation of the MCMC algorithm we adopt for density
estimation and for inferring on the clustering structure of the two
samples. We will be able to skip the evaluation of the sum by resorting
to suitable auxiliary variables whose full conditionals can be
determined and evaluated. To clarify this point, consider the first
sample $\bm{X}^{(n_1)}$, fix $\bm{i}\in\{0,1\}^{k_1}$ and denote by
$\bm
{n}_{(\bm{i})}^0$ the vector of cluster frequencies that correspond to
labels in $\bm{i}$ equal to $0$ whereas $\bm{n}_{(\bm{i})}^1$ is the
vector of cluster frequencies corresponding to labels in $\bm{i}$ equal
to $1$. In a similar fashion, for the second sample $\mathbf{Y}^{(n_2)}$,
for $\bm{l}\in\{0,1\}^{k_2}$, set $\bm{n}_{(\bm{l})}^0$ and $\bm
{n}_{(\bm{l})}^1$. Finally, let $\bm{n}_{(\bm{i},\bm{l})}=(\bm
{n}_{(\bm
{i})}^1,\bm{n}_{(\bm{l})}^1,q_{1,1}+q_{1,2}, \ldots
,q_{k_0,1}+q_{k_0,2})$. From these definitions, it is obvious that $\bm
{n}_{(\bm{i})}^0$, $\bm{n}_{(\bm{l})}^0$ and $\bm{n}_{(\bm{i},\bm{l})}$
are vectors with $k_1-|\bm{i}|$, $k_2-|\bm{l}|$ and $k_0+|\bm
{i}|+|\bm
{l}|$ coordinates, respectively. Moreover, let $\lambda_1$, $\lambda_2$
and $\lambda_0$ be permutations of the coordinates of the vectors $\bm
{n}_{(\bm{i})}^0$, $\bm{n}_{(\bm{l})}^0$ and $\bm{n}_{(\bm{i},\bm
{l})}$. We shall further denote
\[
\Pi_{k,\bm{i},\bm{l}}^{(n_1+n_2)} \bigl(\bm{n}_{(\bm{i})}^0,
\bm{n}_{(\bm{l})}^0,\bm{n}_{(\bm{i},\bm{l})}\bigr)
\]
as the pEPPF conditional on independent random variables $\bm{i}$ and
$\bm{l}$ whose distribution is Bernoulli with parameter $(1-z)$.
Moreover, note that the pEPPF $\Pi_{k}^{(n_1+n_2)}$ depends on the
vectors $\bm{q}^{(i)}$, for $i=1,2$, through their componentwise sum
$\bm{q}^*=(q_{1,1}+q_{1,2}, \ldots,q_{k_0,1}+q_{k_0,2})$. Hence, we
can also write
\[
\Pi_{k}^{(n_1+n_2)}\bigl(\bm{n}^{(1)},
\bm{n}^{(2)},\bm{q}^{(1)},\bm{q}^{(2)}\bigr) =
\Pi_{k}^{(n_1+n_2)}\bigl(\bm{n}^{(1)},\bm{n}^{(2)},
\bm{q}^*\bigr)
\]
%
and shall denote as $\lambda_{1}'$, $\lambda_{2}'$ and
$\lambda_{0}'$ permutations of the components in
$\mathbf{n}^{(1)}$, $\mathbf{n}^{(2)}$ and $\mathbf{q}^*$, respectively. Similarly,
$\lambda_1$, $\lambda_2$ and $\lambda_0$ are permutations of the components
in\vspace*{1.3pt} $\mathbf{n}^{0}_{(\mathbf{i})}$, $\mathbf{n}^{0}_{(\mathbf{l})}$ and $\mathbf{n}_{(\mathbf{i},\mathbf{l})}$. Therefore,
 as a straightforward consequence
of Proposition~\ref{pr2} we obtain the following invariance property for $\Pi
_{k}^{(n_1+n_2)}$ and for $\Pi_{k,\bm{i},\bm{l}}^{(n_1+n_2)}$ whose
proof is omitted since it is immediate.\vspace*{-1pt}

\begin{proposition}\label{pr3}
Let $(\tilde p_1, \tilde p_2)$ be a
GM-dependent normalized CRM defined in \eqref{eqvectornrmi}. Then\vspace*{-1pt}
%
%
\begin{eqnarray}
\label{eqinvarnoindic} \Pi_{k}^{(n_1+n_2)}\bigl(
\bm{n}^{(1)},\bm{n}^{(2)},\bm{q}^*\bigr) &=&
\Pi_{k}^{(n_1+n_2)}\bigl(\lambda_1'
\bm{n}^{(1)},\lambda_2' \bm {n}^{(2)},
\lambda_0' \bm{q}^*\bigr),
\\[-1pt]
\label{eqinvarindic} \Pi_{k,\bm{i},\bm{l}}^{(n_1+n_2)}\bigl(
\bm{n}_{(\bm{i})}^0,\bm {n}_{(\bm{l})}^0,
\bm{n}_{(\bm{i},\bm{l})}\bigr) &=& \Pi_{k,\bm{i},\bm{l}}^{(n_1+n_2)}\bigl(
\lambda_1\bm{n}_{(\bm{i})}^0, \lambda_2
\bm{n}_{(\bm{l})}^0,\lambda_0 \bm{n}_{(\bm{i},\bm{l})}
\bigr).
\end{eqnarray}
\end{proposition}

The invariance property in \eqref{eqinvarnoindic} entails that
exchangeability holds true within three separate groups of clusters:
those with nonshared values and the clusters shared by the two
samples. Such a finding is not a surprise since it reflects the partial
exchangeability assumption.
On the other hand, \eqref{eqinvarindic} implies that, conditional on
a realization of $\bm{i}$ and $\bm{l}$ whose components are i.i.d.
Bernoulli random variables with parameter $1-z$, a similar partially
exchangeable structure is revealed even if it now involves different
groupings of the clusters that are still three: two groups with
nonshared values that are labeled either by $i_j$ or $l_j$ equal to
$0$, and the group containing both observations shared by the two
samples and nonshared values labeled by either $i_j$ or $l_j$ equal to
$1$. 
Moreover, unlike \eqref{eqinvarnoindic} these three groups of
clusters are governed by independent random probability measures.
The invariance structure displayed in~\eqref{eqinvarindic}
corresponds to a mixture decomposition for $\pt_1$ and $\pt_2$ that is
going to be displayed in the next section and is also relevant in
simplifying the MCMC sampling scheme we are going to devise. Note that
\eqref{eqinvarnoindic} holds true since the sum appearing in the
representation of $\Pi_k^{(n_1+n_2)}$ is over all possible $\{0,1\}
$-valued indices $i_j$ and $l_j$: hence a permutation of the frequency
vectors within the three groups simply yields a permutation of the
summands in Proposition~\ref{pr2}. On the contrary, fixing the indices $i_j$
and $l_j$ as in~\eqref{eqinvarindic} corresponds to dropping the sum
in $\Pi_k^{(n_1+n_2)}$ and, then, the invariance is restricted to those
frequencies that correspond to the same index values.\vspace*{-1pt}

\setcounter{example}{0}
\begin{example}[(Continued)]
Let $(\mt_1,\mt_2)$ be a vector of
GM-dependent gamma CRMs. If $\bm{i}=(i_1,\ldots,i_{k_1})\in\{0,1\}
^{k_1}$ and $\bm{l}=(l_1,\ldots,l_{k_2})\in\{0,1\}^{k_2}$ define
$\bar
{n}_1=\sum_{j=1}^{k_1}(1-i_j)n_{j,1}$, $\bar{n}_2=\sum_{j=1}^{k_2}(1-l_j)n_{j,2}$, $\bar{n}_{1,0}=\sum_{j=1}^{k_1}i_j
n_{j,1}$. Moreover, to further simplify notation, set\vspace*{-2pt}
\[
\xi_\sigma\bigl(\bm{n}^{(1)},\bm{n}^{(2)},\bm{q}^*
\bigr)=\prod_{j=1}^{k_1}(1-
\sigma)_{n_{j,1}-1} \prod_{i=1}^{k_2} (1-
\sigma)_{n_{i,2}-1} \prod_{r=1}^m(1-
\sigma)_{q_{r,1}+q_{r,2}-1},
\]
$\alpha'=c+cz+|\bm{q}^*|$ and $\beta'=c+\bar{n}_{1,0}+|\bm{q}^*|$.
%
%
It can then be shown that the pEPPF of the GM-dependent Dirichlet
process is then given by\vspace*{-1pt}
\begin{eqnarray*}
&&\Pi_{k}^{(n_1+n_2)}\bigl(\bm{n}^{(1)},
\bm{n}^{(2)},\bm{q}^*\bigr)
\\[-1pt]
&&\quad= c^k \xi_0\bigl(
\bm{n}^{(1)},\bm{n}^{(2)},\bm{q}^*\bigr)
\\[-1pt]
&&\qquad{}\times \sum_{(*)} 
\frac{z^{k_1+k_2-|\bm{i}|-|\bm{j}|}(1-z)^{k_0+|\bm{i}|+|\bm{j}|}}{
(\alpha')_{n_1}
(\beta')_{n_2}}
{}_3F_2\bigl(cz+\bar{n}_2,
\beta',n_1;n_1+\alpha',n_2+
\beta';1\bigr)
\end{eqnarray*}
for any $(\bm{n}^{(i)},\bm{q}^{(i)})\in\Delta_{n_i,k_i,k_0}$, for
$i=1,2$, and for any $k_1\le n_1$, $k_2\le n_2$ and $k_0$ such that
$k=k_1+k_2+k_0\in\{1,\ldots,n_1+n_2\}$. Note also that if there is only
one sample, namely $n_1n_2=0$, the previous pEPPF reduces to the EPPF
of the Dirichlet process determined in \cite{ewens,antoniak}.
\end{example}

\begin{example}[(Continued)]
When $(\mt_1,\mt_2)$ is a vector of
GM-dependent $\sigma$-stable CRMs, one obtains a pEPPF of the form
\begin{eqnarray*}
&&\Pi_{k}^{(n_1+n_2)}\bigl(\bm{n}^{(1)},
\bm{n}^{(2)},\bm{q}^*\bigr)\\
&&\!\!\quad= \frac
{\Gamma
(k)}{\Gamma(n_1)\Gamma(n_2)} \sigma^{k-1}
\xi_\sigma\bigl(\bm{n}^{(1)},\bm{n}^{(2)},\bm{q}^*
\bigr)
\\
&&\!\!\qquad{}\times \sum_{(*)}z^{k_1+k_2-|\bm{i}|-|\bm{l}|}
(1-z)^{k_0+|\bm{i}|+|\bm
{l}|}%
\int_0^1
\frac{w^{n_1-\bar{n}_{1}+(k_1-|\bm{i}|)\sigma-1}(1-w)^{n_2-\bar
{n}_2+(k_2-|\bm{l}|)\sigma-1}} {
[1-z+zw^\sigma+z(1-w)^\sigma]^k} \,\ddr w,
\end{eqnarray*}
where $\bar{n}_1$, $\bar{n}_2$, $\bar{n}_{1,0}$ are defined as in
Example~\ref{ex1}. Note that the one-dimensional integral above has the same
structure as the one appearing in $\mathcal{I}(c,z)$ and can be
evaluated numerically. Also in this case, if $n_1n_2=0$ the above
expression reduces to the EPPF of the normalized $\sigma$-stable
process. See, for example, \cite{pitman06}.
\end{example}

\begin{remark}\label{rem1}
Following a request of the referees, we also sketch
the extension to more than a pair of dependent random probability
measures the most natural being $\tilde\mu_i=\mu_i+\mu_0$, for each
$i=1,\ldots,N$ and $N>2$. If the mutually independent CRMs $\mu_i$ are
identical in distribution, for $i=1,\ldots,N$, and independent from the
common source of randomness $\mu_0$, one immediately obtains that the
joint Laplace transform of the vector $(\mt_1,\ldots,\mt_N)$ evaluated
at a vector function $(f_1,\ldots,f_N)$ is given by
\[
\E \bigl[\edr^{-\sum_{i=1}^N \mt_i(f_i)} \bigr]=\edr^{-c(1-z)
\psi(|\bm
{f}|)-cz\sum_{i=1}^N \psi(f_i)},
\]
where $\psi$ is the Laplace exponent defined in \eqref{eqlapl} and
shared by the $\mu_i$'s ($i=0,1,\ldots,N$) and $|\bm{f}|=\sum_{i=1}^N
f_i$.
This expression can be used to mimic the proof of Proposition~\ref{pr2}
 and leads
to a straightforward generalization of the pEPPF in the $N$-dimensional
case, which turns out to have the following form
%
\begin{eqnarray*}
&&\frac{c^k}{\prod_{j=1}^N \Gamma(n_j)} \sum_{(*)}(1-z)^{k_0+\sum
_{j=1}^N |\bm{i}_j|}
z^{\sum_{j=1}^N(k_j-|\bm{i}_j|)} \\
&&\quad{}\times\int_{(0,\infty)^N} \prod
_{j=1}^N u_j^{n_j-1}
 \edr^{-c(1-z) \psi(|\bm{u}|)-cz\sum_{i=1}^N\psi(u_i)} \prod_{j=1}^N
\prod_{l=1}^{k_j}\tau_{i_{l,j}}(u_j)
\tau_{n_{l,j}-i_{l,j}}\bigl(|\bm {u}|\bigr)\\
&&\hspace*{32pt}\qquad{}\times\prod_{l=1}^{k_0}
\tau_{q_l}\bigl(|\bm{u}|\bigr) \,\ddr u_1,\ldots,\,\ddr
u_N,
\end{eqnarray*}
where the $(*)$ is the set of all vectors $\bm{i}_j=(i_{1,j},\ldots
,i_{k_j,j})\in\bigtimes_{l=1}^{k_j}\{0,n_{l,j}\}$, for $j=1,\ldots,N$,
$|\bm{u}|=\sum_{i=1}^N u_i$ and $|\bm{i}_j|=\sum_{l=1}^{k_j}i_{l,j}$.
Moreover, the definition of $\tau_q$ in \eqref{eqtau} is extended to
cover the case with $q=0$ as $\tau_0(u)= 1$
for any $u>0$. The previous expression provides the probability of
observing an array of $N$ samples, with respective sizes $n_1,\ldots
,n_N$, with observations partitioned into $k_j$ clusters specific to
the $j$th sample and $k_0$ groups shared by two or more samples. The
exact evaluation of the above $N$-dimensional integral poses some
additional challenges and its implementation within a sampling scheme
is more demanding. A notable exception is given by the GM-dependent
Dirichlet process where for computational purposes one can avoid the
use of the pEPPF and rely on a mixture representation of $\pt_i$ and
$\pt_0$ that will be detailed at the beginning of the next section.
\end{remark}

\section{Dependent mixtures}\label{secddp}
We now apply the general results for GM-dependent normalized CRMs
$(\pt
_1,\pt_2)$ to mixture models with random dependent densities. In fact,
we consider data that are generated from random densities $\tilde f_1$
and $\tilde f_2$ defined by $\tilde f_i(x)
=\int_\Theta h_i(x;\theta) \tilde p_i(\ddr\theta)$, for $i=1,2$, with
$\Theta$ being a complete and separable metric space equipped with the
corresponding Borel $\sigma$-algebra. If $\bbm{\theta}^{(i)}=(\theta
_{1,i},\ldots,\theta_{n_i,i})$, for $i=1,2$, stand for vectors of
latent variables corresponding to the two samples, the mixture model
can be represented in hierarchical form as
%
%
%
\begin{eqnarray}\label{eqgerobserv}
X_{i,1} | \bigl(\bbm{\theta}^{(1)},\bbm{\theta}^{(2)}
\bigr) &\simind& h_1( \cdot;\theta_{i,1}),\qquad i=1,
\ldots,n_1,
\nonumber\\
Y_{j,2} | \bigl(\bbm{\theta}^{(1)},\bbm{\theta}^{(2)}
\bigr) &\simind& h_2( \cdot;\theta_{j,2}),\qquad j=1,
\ldots,n_2,
\nonumber\\
(X_{i,1},Y_{j,2}) | \bigl(\bbm{\theta}^{(1)},\bbm{
\theta}^{(2)}\bigr) &\simind &h_1( \cdot;\theta_{i,1})
h_2( \cdot;\theta_{j,2}),
\\
\theta_{j,i} | (\tilde p_1,\tilde p_2) &
\simiid&\tilde p_i,\qquad j=1,\ldots,n_i; i=1,2,
\nonumber\\
(\tilde p_1,\tilde p_2) &\stackrel{d} {=}&\mbox{GM-dependent normalized
CRM}.\nonumber
\end{eqnarray}
Henceforth, we will set $h_1=h_2=h$; the case of $h_1\ne h_2$ can be
handled in a similar fashion, with the obvious variants. The
investigation of distributional properties of the model is eased by
rewriting $\tilde p_1$ and $\pt_2$ in the following mixture form
%
%
\begin{equation}
\tilde p_i=w_i p_i+(1-w_i)
p_0,\qquad i=1,2, \label{eqmixtrpm}
\end{equation}
where $w_i=\mu_i(\X)\{\mu_i(\X)+\mu_0(\X)\}^{-1}$, 
the $p_i$'s and $p_0$ are independent normalized CRMs with L\'evy
intensities $cz P_0(\ddr x) \rho(s)\,\ddr s$ and $c(1-z) P_0 (\mathrm{d}x)\rho
(s)\,\ddr s$, respectively. Obviously $w_1$ and $w_2$ are dependent. In
general, the weights $w_i$ and the $p_i$'s are dependent, the only
exception being the case in Example~\ref{ex1} where the $p_i$'s are independent
Dirichlet processes. Details about this special case will be provided later.

\begin{remark}\label{rem2}
An interesting aspect of \eqref{eqmixtrpm} is that each $\tilde p_i$
can be decomposed into two independent sources of randomness: an
idiosyncratic one, $p_i$, and a common one, $p_0$. This is close in
spirit to the model of M\"uller, Quintana and Rosner \cite{mqr04},
which is based on a vector of dependent random probability measures
$(\tilde p_1, \ldots, \tilde p_n)$ defined as
%
%
\begin{equation}
\label{eqmueller} \tilde p_i= \omega p_i + (1-\omega)
p_0,
\end{equation}
where $p_i$ and $p_0$ are independent Dirichlet processes and the
distribution of $\omega$ is a mixture with point masses $\omega=0$ and
$\omega=1$ and the remaining mass spread on $(0,1)$ through a beta
density. Despite their similarity, there are however some crucial
differences among GM-dependent normalized CRMs and the model in \eqref
{eqmueller} so that it is not possible to interpret one as the
generalization of the other, nor viceversa. The first thing to note is
that \eqref{eqmueller} assumes common weights, $\omega$ and
$1-\omega
$, for each $\tilde p_i$ whereas in our proposal the weights of the
mixtures $w_i$ in \eqref{eqmixtrpm} do not coincide for different $i$
even if they have the same marginal distributions. More importantly,
the random probability measures defined in \cite{mqr04} via \eqref
{eqmueller} are, in general, marginally not Dirichlet processes. In
our framework, preserving the marginal Dirichlet structure or, in
general, a normalized CRM structure is relevant: it guarantees the
degree of analytical tractability we need for determining
distributional results and devising suitable sampling strategies. The
latter can then be thought of as alternative to the existing algorithms
for dependent random probability measures such as, for example, the one
proposed in \cite{mqr04}.
\end{remark}

On the basis of the decomposition displayed in \eqref{eqmixtrpm}, one
can introduce two collections of auxiliary random variables, $(\zeta
_{j,1})_{j\ge1}$ and $(\zeta_{j,2})_{j\ge1}$, defined on $(\Omega
,\Fcr
,\P)$ and taking values in $\{0,1\}^\infty$ and $\{0,2\}^\infty$, and
provide an useful alternative representation of the mixing measure in~\eqref{eqgerobserv} in terms of these auxiliary variables as
%
%
\begin{eqnarray}
\theta_{i,1} | \zeta_{i,1},\mu_1,
\mu_2,\mu_0 & \simind& p_{\zeta_{i,1}},\qquad i=1,
\ldots,n_1,
\nonumber
\\
\theta_{j,2} | \zeta_{j,2},\mu_1,
\mu_2,\mu_0 & \simind& p_{\zeta_{j,2}},\qquad j=1,
\ldots,n_2, \label{eqgerarchico}
\\
(\zeta_{i,1},\zeta_{j,2}) | \mu_1,
\mu_2,\mu_0 & \simind&\operatorname{bern}
\bigl(w_1;\{0,1\} \bigr) \times\operatorname{bern}
\bigl(w_2;\{0,2\} \bigr),
\nonumber
\end{eqnarray}
where $X\sim\operatorname{bern}(w;\{a,b\})$ means that $\P[X=b]=1-\P[X=a]=w$
for $w\in[0,1]$ and \mbox{$a,b\in\R$}. 
The latent variables $\bbm{\theta}^{(i)}$ are, then, governed by
GM-dependent normalized CRMs. Therefore, we can resort to results
established in Section~\ref{secpepf} to obtain the full conditional
distributions for all the quantities that need to be sampled in order
to attain posterior inferences. Given the structure of the model, the
latent $\bbm{\theta}^{(i)}$, $i=1,2$, might feature ties which generate,
according to the notation we have already introduced, $k_1+k_2+k_0$
clusters. Our analysis of the partition of the $\bbm{\theta}^{(i)}$'s
will further benefit from the following fact that is a straightforward
consequence of Proposition~\ref{pr2}.

\begin{coro}\label{coro1}
Let $(\tilde p_1, \tilde p_2)$
be a GM-dependent normalized CRM defined in \eqref{eqvectornrmi}.
Suppose $P_0$ in (\ref{eqintensity}) is a nonatomic probability measure on $(\X,\Xf)$. Then
%
%
\begin{equation}
\label{eqlabel} \P[\theta_{i,1}=\theta_{j,2} |
\zeta_{i,1}\ne\zeta_{j,2}]=0.
\end{equation}
\end{coro}

Hence, \eqref{eqlabel} entails that ties between the two groups $\bbm
{\theta}^{(1)}$ and $\bbm{\theta}^{(2)}$ may arise with positive
probability only if any two $\theta_{i,1}$ and $\theta_{j,2}$ share the
same label $\zeta_{i,1}=\zeta_{j,2}=0$. This is a structural property
of the model and it intuitively means that there cannot be overlaps
between the different sources of randomness involved, which seems desirable.

Suppose $\bbm{\theta}_*^{(i)}=(\theta_{1,i}^*,\ldots,\theta_{k_i,i}^*)$,
for $i=1,2$, and $\bbm{\theta}_*=(\theta_1^*,\ldots,\theta_{k_0}^*)$
denote the vectors of unique distinct values associated to the
$K=K_1+K_2+K_0$ clusters. The corresponding partition is
%
%
\begin{equation}
\wt{\bbm{\pi}}_{n_1,n_2}=\bigcup_{i=1}^2
\{C_{j,i}\dvt j=1,\ldots,K_i\} \cup \{C_{j,i,0}\dvt
j=1,\ldots,K_0\}, \label{eqpartition}
\end{equation}
where $r,s\in C_{j,i}$ means that $\theta_{r,i}=\theta_{s,i}=\theta
_{j,i}^*$, whereas $r_1\in C_{j,1,0}$ and $r_2\in C_{j,2,0}$ implies
that $\theta_{r_1,1}=\theta_{r_2,2}=\theta_j^*$. It is clear, from the
specification of the model \eqref{eqgerobserv}, that the conditional
density of the data $(\bm{X}^{(n_1)},\bm{Y}^{(n_2)})$, given the
partition $\wt{\bbm{\pi}}_{n_1.n_2}=\bbm{\pi}_{n_1,n_2}$ and the distinct
latent variables $\bbm{\theta}^*=(\bbm{\theta}_*^{(1)},\bbm{\theta
}_*^{(2)},\bbm{\theta}_*)$, coincides with
%
%
\begin{eqnarray}
\label{eqbigprod}&& f\bigl(\bm{x},\bm{y} | \bbm{\theta}^*,\bbm{
\pi}_{n_1,n_2}\bigr)
\nonumber
\\[-8pt]
\\[-8pt]
\nonumber
&&\quad= \prod_{j=1}^{k_1}
\prod_{i\in C_{j,1}}h\bigl(x_i;
\theta_{j,1}^*\bigr)
\prod_{\ell=1}^{k_2} \prod
_{i\in C_{\ell,2}}h\bigl(y_i;\theta_{\ell,2}^*\bigr)
\prod_{r=1}^{k_0} \prod
_{i\in C_{r,1,0}}h\bigl(x_i;\theta_r^*\bigr)
\prod_{\ell\in
C_{r,2,0}}h\bigl(y_\ell;
\theta_r^*\bigr).
\end{eqnarray}
Finally, set
%
%
\begin{equation}
\label{eqlstar} \mathcal{L}^*(\ddr\bm{x},\ddr\bm{y},\ddr\pi,\ddr\bbm{\theta },
\ddr\bbm{\zeta})
\end{equation}
as the distribution of the data $(\bm{X}^{(n_1)}, \bm{Y}^{(n_2)})$, the
partition $\wt{\bbm{\pi}}_{n_1,n_2}$ in \eqref{eqpartition}, the vector
of unique values in $\bbm{\theta}=(\bbm{\theta}^{(1)},\bbm{\theta}^{(2)})$
and the labels $\bbm{\zeta}=(\bbm{\zeta}^{(1)},\bbm{\zeta}^{(2)})$. If
$n=n_1+n_2$, then $\mathcal{L}^*$ is a probability distribution on the
product space $\X^n\times\mathcal{P}_{n}\times\Theta^n\times\{
0,1\}
^{n_1}\times\{0,2\}^{n_2}$, where $\mathcal{P}_{n}$ is the space of all
possible realizations of the random partition $\wt{\bbm{\pi}}_{n_1,n_2}$
in \eqref{eqpartition}. The determination of $\mathcal{L}^*$ will be
first given for any pair of GM-dependent normalized CRMs. The specific
expressions valid for dependent mixtures of the Dirichlet and the
normalized $\sigma$-stable processes will be established as
straightforward corollaries. In the sequel, we also denote as $g_0$ a
density of $P_0$ with respect to some $\sigma$-finite dominating
measure $H$ on $\Theta$, namely $g_0=\ddr P_0/\ddr H$.

\begin{proposition}\label{pr4}
Let $(\tilde p_1, \tilde p_2)$ be a
GM-dependent normalized CRM defined in \eqref{eqvectornrmi}.
Moreover, let $\bbm{\zeta}_i^*=(\zeta_{1,i}^*,\ldots,\zeta_{k_i,i}^*)$
be the vectors of labels corresponding to the distinct latent variables
$\bbm{\theta}_*^{(i)}$, with $i=1,2$. For the dependent mixture model in
\eqref{eqgerobserv}, the distribution $\mathcal{L}^*$ in \eqref
{eqlstar} has density given by
%
%
\begin{equation}
\label{eqjointall} g\bigl(\bm n^{(1)},\bm n^{(2)},\bm
q^{(1)},\bm q^{(2)},\bbm{\zeta}^*\bigr) f\bigl(\bm {x},\bm{y} |
\bbm{\theta}^*, \bbm{\pi}_{n_1,n_2}\bigr) \prod
_{i=1}^k g_0\bigl(\theta_i^*
\bigr),\
\end{equation}
where
%
%
\begin{eqnarray}
\label{eqconditeppf}&& g\bigl(\bm n^{(1)},\bm n^{(2)},\bm
q^{(1)},\bm q^{(2)},\bbm{\zeta}^*\bigr)\nonumber\\
&&\quad= \frac{c^{k} z^{\tilde k_1+\tilde k_2}
(1-z)^{k_1+k_2-\tilde k_1-\tilde k_2}}{\Gamma(n_1)\Gamma(n_2)}
\nonumber\\
&&\qquad{}\times\int_0^\infty\int_0^\infty
u^{n_1-1} v^{n_2-1}\edr^{-cz[\psi
(u)+\psi(v)]
-c(1-z)\psi(u+v)}
\\
&&\hspace*{48pt}\qquad{}\times\prod_{j=1}^{k_1}\tau_{n_{j,1}}
\bigl(u+\bigl(1-\zeta_{j,1}^*\bigr) v\bigr) \prod
_{j=1}^{k_2}\tau_{n_{j,2}}\bigl(\bigl(1-
\zeta_{j,2}^*/2\bigr) u+ v\bigr)
\nonumber\\
&&\hspace*{48pt}\qquad{}\times \prod_{r=1}^{k_0}
\tau_{q_{r,1}+q_{r,2}}(u+v) \,\ddr u \,\ddr v,\nonumber
\end{eqnarray}
where $\tilde k_1=|\bbm{\zeta}_*^{(1)}|$ and $\tilde k_2=|\bbm{\zeta
}_*^{(2)}|/2$ identify the number of clusters with label $1$ and $2$,
respectively.
\end{proposition}

Before examining the details of the models we will refer to for
illustrative purposes, it should be recalled that our approach yields
posterior estimates of $\tilde f_1$ and $\tilde f_2$ and of the number
of clusters $K_{\bm{X}}$ and $K_{\bm{Y}}$ into which one can group the
two sample data. Another interesting issue concerns the estimation of
statistical functionals of $\tilde f_1$ and of $\tilde f_2$, which has
been addressed in the exchangeable case by Gelfand and Kottas \cite
{gk02}. Their approach is based on a suitable truncation of the
stick-breaking representation of the Dirichlet process. In order to
extend their techniques to this setting, a representation of the
posterior distribution of a pair of GM-dependent normalized CRMs is
still missing.

\subsection{Dependent mixtures of Dirichlet processes}\label{sec4.1}
If the vector $(\pt_1, \pt_2)$ is a GM-dependent Dirichlet process as
in Example~\ref{ex1}, then one finds out that the weights $(w_1,w_2)$ in \eqref
{eqmixtrpm} and the Dirichlet process components $p_i$, for
$i=0,1,2$, are independent and the density function of the vector
$(w_1,w_2)$ is
%
%
\begin{equation}
f(w_1,w_2)=\frac{\Gamma(c+cz)}{\Gamma^2(cz)\Gamma(c-cz)} \frac
{(w_1w_2)^{cz-1}
[(1-w_1)(1-w_2)]^{c-1}}{(1-w_1w_2)^{c+cz}}
\indic_{[0,1]^2}(w_1,w_2). \label{eqolkinliu}
\end{equation}
This corresponds to the bivariate beta distribution introduced \cite
{ol}. This model is analyzed in~\cite{cipriota,raoteh}, where
independence between $(w_1,w_2)$ and $(p_0,p_1,p_2)$ is used to devise
a sampler that includes sampling the weights $w_i$. Here we marginalize
with respect to both the weights $(w_1,w_2)$ and the random independent
Dirichlet processes $p_i$, for $i=0,1,2$. The first marginalization is
trickier and is achieved by virtue of the results in Section~\ref{secpepf}.

\begin{coro} \label{prpjointlab}
Let $(\tilde p_1, \tilde p_2)$ be a GM-dependent Dirichlet
process. A density of the probability distribution $\mathcal{L}^*$
defined in \eqref{eqlstar} coincides with
\begin{eqnarray*}
&&\frac{c^kz^{\tilde k_1+\tilde k_2}(1-z)^{k_1+k_2-\tilde k_1-\tilde k_2}} {
(\alpha)_{n_1}(\beta)_{n_2}} {}_3F_2(\alpha-cz+n_1-
\bar{n}_1,n_1,n_2; \alpha+n_1,
\beta+n_2;1)\\
&&\quad{}\times\xi_0\bigl(\bm{n}^{(1)},\bm{n}^{(2)},
\bm{q}^*\bigr) f\bigl(\bm{x},\bm{y} | \bbm{\theta}^*,\bbm{\pi}_{n_1,n_2}
\bigr) 
\prod_{i=1}^k
g_0\bigl(\theta_i^*\bigr),
\end{eqnarray*}
where $\bar{n}_1=|\bbm{\zeta}^{(1)}|$, $\bar{n}_2=|\bbm{\zeta}^{(2)}|/2$,
$\alpha=c+n_2-\bar{n}_2$ and $\beta=c+n_1-\bar{n}_1$.
\end{coro}

As for the actual implementation of the model, a Gibbs sampler easily
follows from Corollary~\ref{prpjointlab}. A key issue is the sampling of the labels.
This can be done by first observing the following facts: (i) if $\theta
_{i,1}=\theta_{j,2}$ then, by Corollary~\ref{coro1}, the corresponding labels are
zero, namely $\zeta_{i,1}=\zeta_{j,2}=0$; (ii) given the partition
$\bbm{\pi}$, the dimensions of label vectors can be shrunk so that one
basically has $k$ labels corresponding to the $k=k_1+k_2+k_0$ clusters
of the partition. Remark (i) implies that we do not need to sample the
labels associated to $\theta_{i,1}$ values coinciding with any of the
$\theta_{j,2}$'s and viceversa. Moreover, remark (ii) implies that for
any $r,s\in C_{j,i}$ one has $\zeta_{r,i}=\zeta_{s,i}$ and, thus, we
need to sample only labels $\zeta_{j,i}^*$ corresponding to distinct
values $\theta_{j,i}^*$. Finally, there might be $\theta_{j,1}$'s (or
$\theta_{j,2}$'s) associated to $\zeta_{j,1}=0$ (or $\zeta_{j,2}=0$)
that do not coincide with any of the $\theta_{i,2}$'s (or of the
$\theta
_{i,1}$'s): the corresponding labels are not degenerate and must be
sampled from their full conditionals. If $\bbm{\zeta}_{-j,*}^{(1)}$
stands for the vector $\bbm{\zeta}_*^{(1)}$ with the $j$th component
removed, we use the short notation
\[
\pi_{j,1}(x):= \P\bigl[\zeta_{j,1}^*=x | \bbm{
\zeta}_{-j,*}^{(1)},\bbm{\zeta }_*^{(2)},\bbm{
\theta}^*,\bm{X}^{(n_1)},\bm{Y}^{(n_2)}\bigr].
\]
Hence, if $\theta^*_{j,1}$ does not coincide with any of the distinct
values of the latent variables for the second sample, it can be easily
deduced that
%
%
\begin{eqnarray}
\label{eqfullcondlabel} \pi_{j,1}(x)&\propto&\indic_{ \{0,1 \}}(x)
\frac{z^x (1-z)^{1-x}}{(\alpha)_{n_2}(\beta_x)_{n_2}}
\nonumber
\\[-8pt]
\\[-8pt]
\nonumber
&&{}\times {}_3F_2 (\alpha-cz+n_1-
\bar{n}_{-j,1}-xn_{j,1},n_1,n_2;\alpha
+n_1,\beta_x+n_2;1 ),
\end{eqnarray}
where $\bar{n}_{-j,1}:=\sum_{i\ne j}n_{i,1}\zeta_{i,1}^*$ with
$n_{i,1}$ denoting the size of the cluster identified by $\theta
_{i,1}^*$. Moreover, $\beta_x=c+n_1-\bar{n}_{-j,1}-xn_{j,1}$.
Obviously, the normalizing constant is determined by $\pi_{j,i}(0)+\pi
_{j,i}(1)=1$. The full conditionals for the $\zeta_{j,2}^*$ can be
determined analogously.\vspace*{1pt}

As for the full conditionals of the $\theta_{j,i}$'s, these reduce to
the ones associated to the univariate mixture of the Dirichlet process,
since one is conditioning on the labels $\zeta_{j,i}$ as well. Hence,
one can sample $\theta_{j,1}$ from
%
%
\begin{equation}
\label{eqfullcondittheta} 
w_0 P_{j,1}^*(\ddr\theta)+\sum
_{l\in\mathcal{J}_{-j, \zeta_{j,1}}} w_l\delta_{\tilde\theta_{l,\zeta_{j,1}}}(\ddr\theta),
\end{equation}
where $\tilde\theta_{l,\zeta_{j,1}}$ are the distinct $\theta$ values
in the urn labeled $\zeta_{j,1}$ and $\mathcal{J}_{-j, \zeta_{j,1}}$ is
the set of indices of distinct values from the urn labeled $\zeta
_{j,1}$ after excluding $\theta_{j,1}$. Moreover,
%
%
\begin{eqnarray}
\label{eqwpred} w_0&\propto& c (1-z)^{1-\zeta_{j,1}} z^{\zeta_{j,1}}
\int_{\Theta} h (x_{j};\theta )P_0(\ddr
\theta),
\nonumber
\\[-8pt]
\\[-8pt]
\nonumber
w_l&\propto& n_{l,1}^{(-j)} h(x_{j};
\tilde\theta_{l,\zeta_{j,1}}).
\end{eqnarray}
In the weights above,\vspace*{-1pt} $P_{j,1}^*(\ddr\theta)=h(x_j;\theta)P_0(\ddr
\theta)/
\int_{\Theta}h(x_j;\theta)P_0(\ddr\theta)$ and $n_{l,1}^{(-j)}$ is the
size of the cluster containing $\tilde{\theta}_{l,\zeta_{j,1}}$, after
deleting $\theta_{j,1}$. With obvious modifications, one also obtains
the full conditional for generating $\theta_{j,2}$. This last point
suggests that, conditional on the labels, one needs to run three
independent Blackwell--MacQueen P\'olya urn schemes: two are related to
the idiosyncratic (and independent) components and one is related to
the common component. Given this, the only difficulty in implementing
the algorithm is due to the generalized hypergeometric function $
{}_3F_2(a,b,c;e,f;x)$. Indeed, when such a function is evaluated at
$x=1$, as in our case, the convergence of the series defining it can be
very slow, depending on the magnitude of $e+f-a-b-c>0$: the lower such
a value, the slower the convergence of the series. The efficiency of
the algorithm can, thus, be improved by suitably resorting to
identities that involve generalized hypergeometric functions in order
to obtain equivalent expressions with a larger value of $e+f-a-b-c$. In
particular, in the examples considered here we have been able to
considerably speed up the implementation of the algorithm by applying
an identity that can be found in \cite{bailey}, page~14.

\subsection{Dependent mixtures of normalized \texorpdfstring{$\sigma$-stable}{sigma-stable} processes}\label{sec4.2}
Consider a GM-dependent normalized $\sigma$-stable CRM vector $(\pt
_1, \pt_2)$ as in Example~\ref{ex2}. The corresponding model is somehow more
complicated to deal with, but at the same time it is more
representative of what happens in the general case since the
simplifications typical of the Dirichlet process do not occur.
Specifically, the weights $(w_1,w_2)$ are no longer independent from
the normalized $\sigma$-stable processes $p_i$ in \eqref{eqmixtrpm}.
Moreover, the density of $(w_1,w_2)$ is not available in closed form
for any $\sigma\in(0,1)$, but only for $\sigma=1/2$. Nonetheless, it is
still possible to obtain analytic forms for the full conditionals
allowing to estimate the marginal densities $\tilde f_i$ and to analyze
the clustering structure featured by the two-sample data. Indeed, one
can show the following
corollary.

\begin{coro}
\label{prpjointlabstable}
Let $(\tilde p_1, \tilde p_2)$ be a GM-dependent normalized
$\sigma$-stable CRM. A density of the probability distribution
$\mathcal{L}^*$ defined in \eqref{eqlstar}
coincides with
\begin{eqnarray*}
&&\frac{z^{\tilde k_1+\tilde k_2}(1-z)^{k_1+k_2-\tilde k_1-\tilde
k_2}}{\Gamma(n_1) \Gamma(n_2)} \sigma^{k-1} \Gamma(k) f\bigl(\bm{x},\bm{y} | \bbm{
\theta}^*,\bbm{\pi}_{n_1,n_2}\bigr) \prod_{i=1}^k
g_0\bigl(\theta_i^*\bigr)
\\
&&\quad{}\times \xi_\sigma\bigl(\bm{n}^{(1)},\bm{n}^{(2)},
\bm{q}^*\bigr) \int_0^1\frac{w^{n-\bar{n}_1+\tilde k_1\sigma-1}(1-w)^{n_2-\bar
{n}_2+\tilde k_2\sigma-1}}{
\{1-z+zw^\sigma+z(1-w)^\sigma \}^k} \,\ddr
w,
\end{eqnarray*}
where $\bar{n}_1=|\bbm{\zeta}^{(1)}|$ and $\bar{n}_2=|\bbm{\zeta
}^{(2)}|/2$.
\end{coro}

In a similar fashion to the dependent Dirichlet process case, from
Corollary~\ref{prpjointlabstable} one can deduce the full conditionals for both the labels
$\bbm{\zeta}^{(i)}_*$ and the $\bbm{\theta}_*^{(i)}$. As for the former,
if $\zeta_{j,1}^*$ corresponds to a distinct value $\theta_{j,i}^*$ not
coinciding with any value $\theta_{l,2}$ from the second sample,
then\vspace*{1pt}
%
%
\begin{eqnarray}
\label{eqfullcondlabelstable} \pi_{j,1}(x)&\propto&\indic_{ \{0,1 \}}(x)
z^x(1-z)^{1-x}
\nonumber
\\[-8pt]
\\[-7pt]
\nonumber
&&{}\times \int_0^1\frac{w^{n-\bar{n}_{-j,1}-xn_{j,1}+(\tilde
k_{-j,1}+x)\sigma-1}(1-w)^{n_2-\bar{n}_2+\tilde k_2\sigma-1}}{
\{1-z+zw^\sigma+z(1-w)^\sigma \}^k} \,\ddr w,
\end{eqnarray}
where $\bar{n}_{-j,1}=\sum_{i\ne j}n_{i,1} \zeta_{i,1}^*$ and
$\tilde
k_{-j,1}=|\bbm{\zeta}_{-j,*}^{(1)}|$.

Interestingly, the full conditionals for the latent random variables
are as simple as in the Dirichlet process case. Since we are again
conditioning on the labels $\bbm{\zeta}^{(1)}$, it is apparent that one
just needs to run three independent Blackwell--MacQueen P\'olya urn
schemes. For $\theta_{j,1}$ the full conditional coincides with \eqref
{eqfullcondittheta} with different weights\vspace*{1pt}
%
%
\begin{eqnarray}
\label{eqweightsstable} w_0&\propto& k_{-j,\zeta_{j,1}}
\sigma(1-z)^{1-\zeta_{j,1}} z^{\zeta_{j,1}} \int_\Theta
h(x_j;\theta) P_0(\ddr\theta),
\nonumber
\\[-8pt]
\\[-8pt]
\nonumber
w_l&\propto&\bigl(n_{l,1}^{(-j)}-\sigma\bigr)
h(x_j;\tilde{\theta}_{l,\zeta_{j,1}}),
\end{eqnarray}
where $k_{-j,\zeta_{j,1}}$ above is the number of clusters associated
to $p_{\zeta_{j,1}}$ after excluding $\theta_{j,1}$.\vspace*{1pt}

\section{Full conditional distributions}\label{sec5}
The results in Sections~\ref{secpeppf} and \ref{secddp} form the basis for the concrete
implementation of the model~\eqref{eqgerobserv}
to a
real datasets in the following section. Here we provide a detailed
description of the algorithm set forth in Section~\ref{secddp} for
specific choices of the kernel $h( \cdot; \cdot)$ and of the
random probability measures $\tilde p_1$ and $\tilde p_2$. In
particular, we make the standard assumption of $h( \cdot;M,V)$ being
Gaussian with mean $M$ and variance $V$ and consider GM-dependent
Dirichlet and normalized $\sigma$-stable processes as mixing measures.
As for the specification of the base measures $P_0$ of such mixing
measures (see \eqref{eqintensity}), we propose a natural extension to
the partially exchangeable case of the quite standard specification of
Escobar and West \cite{ew95}, which greatly contributed to popularizing
the mixture of Dirichlet process model. In particular, we take $P_0$ to
be a normal/inverse-Gamma distribution\vspace*{1pt}
\[
P_0(\ddr M,\ddr V)=P_{0,1}(\ddr V) P_{0,2}(\ddr
M | V)
\]
with $P_{0,1}$ being an inverse-Gamma probability distribution with
parameters $(s,S)$ and $P_{0,2}$ is Gaussian with mean $m$ and variance
$\tau V$. Moreover, the corresponding hyperpriors are of the form\vspace*{1pt}
%
%
\begin{eqnarray}
\label{spec2} \tau^{-1} &\sim&\operatorname{Ga}(w/2,W/2),
\nonumber
\\[1pt]
m &\sim&\mathrm{N}(a,A),
\\
z &\sim &U(0,1),
\nonumber\\
c &\sim&\operatorname{Ga}(a_0,b_0)\nonumber
\end{eqnarray}
for some $w>0$, $W>0$, $A>0$, $a_0>0$, $b_0>0$ and real $a$. In the
following, we focus on the two special cases and provide the analytic
expressions for the corresponding full conditional distributions. In
terms of the notation set in Section~\ref{secddp}, the latent
variables now become $\theta_{j,i}=(M_{j,i},V_{j,i})\in\R\times\R^+$,
for any $j=1,\ldots,n_i$ and $i=1,2$. Moreover, $\tilde\theta
_{j,i}=(\tilde M_{j,i},\tilde V_{j,i})$, for $i=0,1,2$, represent the
$j$th distinct value of the latent variables with label $i$. Also
recall that the number of distinct values with label $i$, for $i=1,2$,
is equal to $\tilde k_i$ and set $\tilde k_0=k_1+k_2-\tilde k_1-\tilde k_2$.

\subsection{GM-dependent Dirichlet processes}\label{sec5.1}
Let us first deal with the hierarchical mixture model \eqref
{eqgerobserv} with $(\pt_1,\pt_2)$ a vector of GM-dependent
Dirichlet processes with parameters $(c,z;P_0)$, which we will denote
by GM--$\mathscr{D}(c,z;P_0)$ in the sequel. With this specification
and the auxiliary variable representation of the mixing measure laid
out in \eqref{eqgerarchico}, the weights of the predictive \eqref
{eqfullcondittheta} are similar to those described in \cite{ew95},
the only differences being related to the bivariate structure, which
results in the dependence on $z$ (see \eqref{eqintensity}) and on the
label $\zeta_{j,i}$. These identify the full conditional for the latent
$\theta_{j,i}$.

In order to determine the full conditionals for the other parameters to
be sampled, let $\bm{D}_{-r}$ stand for the set of all
(hyper)parameters of the model but $r$. As for the full conditional for
$z$, one has
\begin{eqnarray*}
\kappa_z\bigl(z | \bm{X}^{(n_1)},\bm{Y}^{(n_2)},
\bm{D}_{-z}\bigr)&\propto& \kappa_z(z)
z^{\tilde k_1+\tilde k_2}(1-z)^{\tilde k_0}
\\
&&{}\times {}_3F_2 (
\alpha-cz+n_1-\bar{n}_1,n_1,n_2;
\alpha+n_1,\beta +n_2;1 ),
\end{eqnarray*}
%
where $\kappa_z$ is the prior distribution of $z$, which in our
specification coincides with the uniform on $(0,1)$. On the other hand,
an expression for the full conditional for $c$ is obtained as follows
\begin{eqnarray*}
\kappa_c\bigl(c | \bm{X}^{(n_1)},\bm{Y}^{(n_2)},
\bm{D}_{-c}\bigr)&\propto&\kappa_c(c)
\frac{c^{k}}{(\alpha)_{n_1}(\beta)_{n_2}}
\\
&&{}\times {}_3F_2 (
\alpha-cz+n_1-\bar{n}_1,n_1,n_2;
\alpha+n_1, \beta +n_2;1 ),
\end{eqnarray*}
where $\kappa_c$ is the prior distribution of $c$ that is supposed
coincide with $\operatorname{Ga}(a_0,b_0)$. Moreover, note that both the
coefficients $\alpha$ and $\beta$ appearing in the generalized
hypergeometric function ${}_3F_2$ above depend on $c$. See Corollary~\ref{prpjointlab}.
Finally, $\tau$ and $m$ are sampled from the following distributions
%
%
\begin{eqnarray}
\label{fullcondtau} \tau| \bigl(\bm{X}^{(n_1)},\bm{Y}^{(n_2)},
\bm{D}_{-\tau}\bigr) &\sim&\operatorname{IG} \biggl(\frac{w+k}{2},
\frac{W+W'}{2} \biggr),
\\
\label{fullcondm} m | \bigl(\bm{X}^{(n_1)},\bm{Y}^{(n_2)},
\bm{D}_{-m}\bigr) &\sim&\mathrm{N}(RT,T),
\end{eqnarray}
where $\operatorname{IG}(a,b)$ denotes the inverse-gamma distribution with density
function $g(s)\propto  s^{-a-1} \edr^{-\beta/s} \indic_{\R^+}(s)$,
$W'=\sum_{i=0}^2
\sum_{l=1}^{\tilde k_i}(\tilde M_{l,i}-m)^2/\tilde V_{l,i}$
and
\begin{eqnarray*}
T&=& \Biggl[\frac{1}{A}+\frac{1}{\tau} \Biggl(\sum
_{i=1}^{\tilde
k_1}\frac
{1}{\tilde V_{i,1}}+ \sum
_{j=1}^{\tilde k_2}\frac{1}{\tilde V_{j,2}}+ \sum
_{r=1}^{\tilde k_0}\frac{1}{\tilde V_{r,0}} \Biggr)
\Biggr]^{-1},
\\
R&=& \Biggl[\frac{a}{A}+\frac{1}{\tau} \Biggl(\sum
_{i=1}^{\tilde k_1}\frac{\tilde M_{i,1}}{\tilde V_{i,1}}+ \sum
_{j=1}^{\tilde k_2}\frac{\tilde M_{j,2}}{\tilde V_{j,2}}+\sum
_{r=1}^{\tilde k_0}\frac{\tilde M_{r,0}}{\tilde V_{r,0}} \Biggr)
\Biggr]^{-1}.
\end{eqnarray*}

\subsection{GM-dependent \texorpdfstring{$\sigma$-stable}{sigma-stable} normalized random measures}\label{sec5.2}
When $(\pt_1,\pt_2)$ is a vector of GM-dependent normalized $\sigma
$-stable processes with parameters $(z,P_0)$ we set the short notation
GM--st$(\sigma,z,P_0)$. The full conditionals are then derived from
Corollary~\ref{prpjointlabstable}. In particular, explicit expressions for the weights in
\eqref{eqweightsstable} can be deduced and the full conditional for
$z$ which coincides with
\begin{eqnarray*}
\kappa_z\bigl(z | \bm{X}^{(n_1)},\bm{Y}^{(n_2)},
\bm{D}_{-z}\bigr)&\propto&
\kappa_z(z) z^{\tilde k_1+\tilde
k_2}(1-z)^{\tilde
k_0}
\\
&&{}\times\int_0^1\frac{w^{n_1-\bar n_1+\tilde k_1\sigma
-1}(1-w)^{n_2-\bar n_2+\tilde k_2\sigma-1}}{ \{1-z+zw^\sigma
+z(1-w)^\sigma \}^k}\,\ddr w,
\end{eqnarray*}
where $\kappa_z$ is, as in Section~\ref{sec5.1}, uniform on $(0,1)$. Moreover,
if a prior on $(0,1)$ is assigned to the parameter $\sigma$, the
corresponding full conditional is given by
\begin{eqnarray*}
\kappa_\sigma\bigl(\sigma| \bm{X}^{(n_1)},\bm{Y}^{(n_2)},
\bm{D}_{-\sigma}\bigr)&\propto&
\kappa_\sigma(\sigma) \sigma^{k-1} \xi _{\sigma}\bigl(
\bm{n}^{(1)},\bm{n}^{(2)},\bm{q}^*\bigr)
\\
&&{}\times \int_0^1\frac{w^{n_1-\bar n_1+\tilde k_1\sigma
-1}(1-w)^{n_2-\bar n_2+\tilde k_2\sigma-1}}{ \{1-z+zw^\sigma
+z(1-w)^\sigma \}^k}
\,\ddr w. 
\end{eqnarray*}
Finally, the full conditionals for $\tau$ and $m$ coincide with those
displayed in \eqref{fullcondtau} and \eqref{fullcondm} since they
depend only on $h$ and $P_0$ and not on the specific vector of random
probabilities $(\pt_1,\pt_2)$ driving the respective dependent mixtures.

\subsection{Accelerated algorithm}\label{sec5.3}
It is well known that univariate P\'olya urn samplers like the one
proposed in \cite{ew95} tend to mix slowly when the probability of
sampling a new value, $w_0$, is much smaller than the probability to
sample an already observed one. When this occurs, the sampler can get
stuck at the current set of distinct values and it may take many
iterations before any new value is generated. Such a concern clearly
extends also to our bivariate P\'olya urn sampler and, in particular,
to \eqref{eqfullcondittheta} and \eqref{eqweightsstable} leading
the algorithm to get stuck in some specific $\{\tilde\theta_{i,l}\dvt l=0,1,2; i=1,\ldots,\tilde k_l\}$. To circumvent this problem, we
resort to the method suggested in \cite{WesMulEsc94} and \cite
{Mac94}: it consists in resampling, at the end of every iteration, the
distinct values $\tilde\theta_{l,i}$ from their conditional
distribution. Since this distribution depends on the choice of $\pt_1$
and $\pt_2$ only through their base measure $P_0$, it is the same for
the Dirichlet and $\sigma$-stable cases. In particular, for every
$i=1,\ldots,\tilde k_1$, the required full conditional density of
$\tilde\theta_{i,1}$ is
%
%
\begin{equation}
\label{eqfullaccel} \La^*\bigl(\tilde\theta_{i,1} |
\bm{X}^{(n_1)},\bm{Y}^{(n_2)},\bm {D}_{-\tilde\theta_{i,1}}\bigr) \propto
g_0(\tilde\theta_{i,1})\prod_{j\in
C_{i,1}}h(x_{j},
\tilde\theta_{i,1}),
\end{equation}
where $\La^*$ is the joint law defined in \eqref{eqlstar}. 
With our specification, the full conditional distribution of $\tilde
\theta_{i,1}$ in \eqref{eqfullaccel} becomes normal/inverse-Gamma with
\begin{eqnarray*}
\tilde V_{i,1}^{-1} &\sim&\operatorname{Ga} \biggl(s+
\frac{n_{i,1}}{2},S+\frac{\sum x_{j}^2}{2}+\frac
{m^2n_{i,1}-\sum x_{j} (2m+\tau\sum x_{j} )}{2(1+\tau
n_{i,1})} \biggr),
\\
\tilde M_{i,1} | \tilde V_{i,1} &\sim&\mathrm{N} \biggl(
\frac{m+\tau\sum x_{j}}{1+\tau
n_{i,1}},\tilde V_{i,1} \frac{\tau}{1+\tau n_{i,1}} \biggr),
\end{eqnarray*}
where $\sum x_{j}$ is a shortened notation for $\sum_{j\in C_{i,1}}x_{j}$.
Analogous expressions, with obvious modifications, hold true for
$\tilde
\theta_{i,2}$ and $\tilde\theta_{i,0}$.

\section{Illustration}\label{sec6}
In this section, we illustrate the inferential performance of the
proposed model on a two-sample dataset and to this end we implement
the Gibbs sampling algorithm devised in the previous section for \eqref
{eqgerobserv}. We shall consider $(\pt_1,\pt_2)$ being either a
GM--$\Dcr(c,z,P_0)$ or a GM--st$(\sigma,z,P_0)$. In terms of
computational efficiency, we note in advance that the algorithm with
the GM--st mixture is remarkably faster than the one associated to the
GM--$\Dcr$ mixture. As already pointed out in the previous sections,
this is due to the need of repeated evaluations of generalized
hypergeometric function ${}_3F_2$ in the GM--$\Dcr$ case. In contrast,
the numerical evaluation of the one-dimensional integral in
Corollary~\ref{prpjointlabstable}, for the GM--st mixture, is straightforward.

We shall analyze the well-known Iris dataset, which contains measures
of $4$ features of $3$ different species of \textit{Iris} flowers:
\textit{Setosa}, \textit{Versicolor} and \textit{Virginica}. For each
of these species $150$ records of sepal length, sepal width, petal
length and petal width of flowers are available. These data are
commonly used in the literature as an illustrative example for
discriminant analysis. Indeed, it has been noted that \textit{Setosa}
is very well separated from the other two species, which partially
overlap. Of the $4$ measured features, here we consider the petal width
expressed in millimeters. A total number of $50$ observations per
species have been recorded. The $150$ observations are, then, used to
form two samples $\bm{X}^{(n_1)}$ and $\bm{Y}^{(n_2)}$ as follows. We
set $n_1=90$ and let the first sample consist of $50$ observations of
\textit{Setosa} and $40$ of \textit{Versicolor}. Correspondingly
$n_2=60$ and includes $50$ observations of \textit{Virginica} and the
remaining $10$ observations of \textit{Versicolor}. The particular
design of the experiment is motivated by the idea that the \textit
{Versicolor} species identifies the shared component between the two
mixtures, thus making our approach for modeling dependence appropriate.
Moreover, on the basis of previous considerations it is expected that
the two species in the first dataset are more clearly separated than
the two species forming the second sample.

%
\begin{figure}

\includegraphics{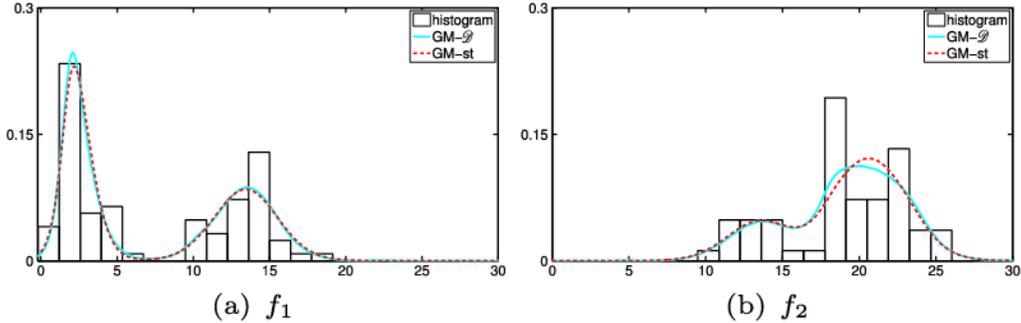}

\caption{GM--$\Dcr(c,z,P_0)$ (solid line) vs.~GM--st$(\sigma,z,P_0)$ (dashed
line) mixture with random $c$ and $\sigma$, respectively: posterior
estimates of the densities $f_1$ and $f_2$. The displayed histograms
are generated by the actual two-sample data.}
\label{irisdensity}
\end{figure}

Our statistical analysis has the following two goals: on the one hand
we wish to estimate the densities generating the two samples and, on
the other, we aim at obtaining an approximation of the posterior
distribution of the number of clusters in each sample. This allows to
draw a direct comparison of the inferential outcomes produced by the
GM--$\Dcr(c,z,P_0)$ and GM--st$(\sigma,z,P_0)$ mixtures. As for the
specifications of the hyperparameters in \eqref{spec2} we essentially
adopted the quite standard specifications of \cite{ew95}. Hence, we
have set $(w,W)=(1,100)$, $(s,S)=(1,1)$, $(a,A)=((n_1\bar{X}+n_2\bar
{Y})/(n_1+n_2), 2)$ and $(a_0,b_0)=(2,1)$ where $\bar{X}$ and $\bar
{Y}$ are the sample means for $\bm{X}^{(n_1)}$ and $\bm{Y}^{(n_2)}$,
respectively. As for the other parameters involved, we suppose that
$c\sim\operatorname{Ga}(2,1)$, whereas $\sigma$ and $z$ are both
uniform on $[0,1]$.
Moreover, these three parameters are independent. All estimates will be
based on $80\mbox{,}000$ iterations of the algorithm after $20\mbox{,}000$ burn-in sweeps.

\begin{figure}[t]

\includegraphics{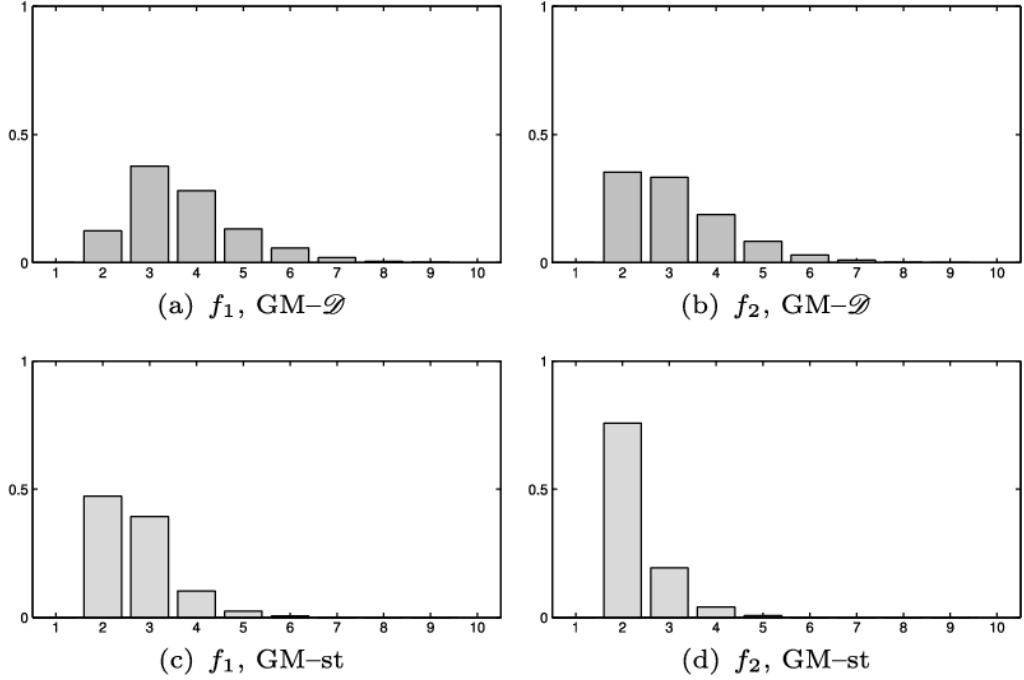}

\caption{GM--$\Dcr(c,z,P_0)$ (top row) vs.~GM--st$(\sigma,z,P_0)$ (bottom row)
mixture with random $c$ and $\sigma$, respectively: posterior
distributions of the number of clusters $K_{\bm{X}}$ and $K_{\bm{Y}}$.}
\label{irisclust}
\end{figure}

The estimated densities are displayed in Figure~\ref{irisdensity} and
there seem to be no significant differences. However, regardless the
particular mixture model specification, the two species forming each
sample are clearly better separated in the first sample. This is not
surprising, given that the second sample is formed by two overlapping
species. See also the histogram in the background of Figure~\ref{irisdensity}.
The results on the clustering structure are reported in Figure~\ref{irisclust} and in Table~\ref{tableIRIS}.
Figure~\ref{irisclust} shows that the posterior distributions of the
number of clusters corresponding to the GM--st mixture is characterized
by a lower variability than in the GM--$\Dcr$ mixture case. Moreover,
if one roughly thinks of each species of flowers in a sample as forming
a single cluster, then it is apparent that the GM--st mixture better
estimates both $K_{\bm{X}}$ and $K_{\bm{Y}}$. See also Table~\ref{tableIRIS}. These results seems to suggest that the parameter $\sigma
$, associated to the stable CRM, has a beneficial impact on the
estimation of the clustering structure. This is in line with the
findings of \cite{lmp} in the exchangeable case, where it is pointed
out that $\sigma$ induces a reinforcement mechanism which improves the
capability of learning the clustering structure from the data. We
believe this aspect is of great relevance and, hence, deserves further
investigation.

%

%
\begin{table}[b]
\caption{GM--$\Dcr(c,z,P_0)$ vs.~GM--st$(\sigma,z,P_0)$ mixture with random $c$
and $\sigma$, respectively: estimated number of clusters (Col.~1 and
2), maximum a posteriori values ($\hat{K}_{\bm{X}}$ and $\hat
{K}_{\bm
{Y}}$) and probability of more than $4$ clusters per sample (Col.~5 and 6)}
\label{tableIRIS}
\begin{tabular*}{\textwidth}{@{\extracolsep{\fill}}lllllll@{}}
\hline
&$\E[K_{\bm{X}}|\ldots]$&$\E[K_{\bm{Y}}|\ldots]$&$\hat K_{\bm
{X}}$&$\hat K_{\bm{Y}}$&$\P(K_{\bm{X}}\geq4)$&$\P(K_{\bm{Y}}\geq
4)$\\
\hline
GM--$\Dcr(c,z,P_0)$&3.72&3.15&3&2&0.50&0.31\\
GM--st$(\sigma,z,P_0)$&2.70&2.30&2&2&0.13&0.05\\
\hline
\end{tabular*}
\end{table}

\begin{appendix}\label{app}
\section*{Appendix}
\subsection{Proof of Proposition \texorpdfstring{\protect\ref{pr1}}{1}}\label{appa1} By combining the definition
of GM-dependent normalized CRMs given in \eqref{eqvectornrmi} with
the gamma integral, it is possible to write
\[
\E \bigl[\tilde p_1(A) \tilde p_2(B) \bigr]=\int
_0^\infty\int_0^\infty
\E \bigl[\edr^{-u\tilde\mu_1(\X)-v\tilde\mu_2(\X)} \tilde\mu_1(A) \tilde
\mu_2(B) \bigr] \,\ddr u\,\ddr v.
\]
Since $\tilde\mu_i=\mu_i+\mu_0$ for $i=1,2$, with $\mu_0$, $\mu
_1$ and
$\mu_2$ independent, one has
\begin{eqnarray*}
&&\E \bigl[\edr^{-u\mu_1(\X)-(u+v)\mu_0(\X)} \mu_1(A) \mu _0(B) \bigr]
\\
&&\quad= \E \bigl[\edr^{-u\mu_1(\X)} \mu_1(A) \bigr] \E \bigl[\edr
^{-(u+v)\mu
_0(\X)} \mu_0(B) \bigr]
\\
&&\quad=c^2 z(1-z)P_0(A) P_0(B)
\edr^{-cz\psi(u)-c(1-z)\psi(u+v)} \tau _1(u) \tau_1(u+v).
\end{eqnarray*}
Use the symbol $A^i$ to denote $A$ if $i=1$ and $A^c$ if $i=0$. Hence,
$\{A^i\cap B^j\dvt i,j=0,1\}$ is the partition of $\X$ generated by $\{
A,B\}$. Hence,
\[
\E \bigl[\edr^{-(u+v)\mu_0(\X)} \mu_0(A) \mu_0(B) \bigr] =
\sum_{i,j=0}^1 \E \bigl[
\edr^{-(u+v)\mu_0(\X)} \mu_0\bigl(A\cap B^i\bigr) \mu
_0\bigl(A^j\cap B\bigr) \bigr].
\]
This implies that
\begin{eqnarray*}
&&\E \bigl[\edr^{-(u+v)\mu_0(\X)} \mu_0(A)\mu_0(B) \bigr] \\
&&\quad=
\edr^{-c(1-z)\psi(u+v)} c(1-z)\\
&&\qquad{}\times \Biggl\{ P_0(A\cap B)
\tau_2(u+v) +c(1-z) \tau_1^2(u+v)\sum
_{i,j=0}^1 P_0\bigl(A\cap
B^i\bigr)P_0\bigl(A^j\cap B\bigr)\Biggr\}
\\
&&\quad=\edr^{-c(1-z)\psi(u+v)} c(1-z) \bigl\{ P_0(A\cap B)
\tau_2(u+v)
+c(1-z)P_0(A)P_0(B)\tau_1^2(u+v)
\bigr\}.
\end{eqnarray*}
Summing up, it follows that
\begin{eqnarray*}
&&\E \bigl[\tilde p_1(A) \tilde p_2(B) \bigr]\\
&&\quad= \int
_0^\infty\int_0^\infty
\edr^{-z(\psi(u)+\psi(v))-c(1-z)\psi(u+v)}
\\
&&\hspace*{36pt}\qquad{}\times c \bigl\{(1-z)P_0(A\cap B) \tau_2(u+v)
+c^2 P_0(A)P_0(B)\\
&&\hspace*{56pt}\qquad{}\times\bigl[(1-z)^2
\tau_1^2(u+v)+z(1-z)\tau_1(u+v) \bigl(\tau_1(u)+\tau_1(v)
\bigr)\\
&&\hspace*{90pt}{}+z^2 \tau_1(u)\tau_1(v)\bigr]\bigr\}
\ddr u\,\ddr v.
\end{eqnarray*}
If in the previous expression one sets $A=B=\X$, then the following
identity holds true
\begin{eqnarray*}
&&c^2 \int_0^\infty\int
_0^\infty\edr^{-cz(\psi(u)+\psi
(v))-c(1-z)\psi
(u+v))}\\
&&\hspace*{47pt}{}\times \bigl[
(1-z)^2 \tau_1^2(u+v)
+z(1-z)\tau_1(u+v) \bigl(\tau_1(u)+\tau_1(v)
\bigr)+z^2\tau_1(u)\tau_1(v)\bigr] \,\ddr u \,
\ddr v
\\
&&\quad= 1-c(1-z) \int_0^\infty\int_0^\infty
\edr^{-c(1-z)\psi(u+v)} \tau _2(u+v) \,\ddr u\,\ddr v.
\end{eqnarray*}
The results in \eqref{eqmixmom} and in \eqref{eqcorrelation} then
follows.

\subsection{Proof of Proposition \texorpdfstring{\protect\ref{pr2}}{2}}\label{appa2} We first determine the
probability distribution of $(\tilde\pi_{n_1,n_2},\bm{X}^{(n_1)},\bm
{Y}^{(n_2)})$. Here $\tilde\pi_{n_1,n_2}$ denotes a random partition of
$\{\bm{X}^{(n_1)},\bm{Y}^{(n_2)}\}$ whose generic realization, $\pi
_{n_1,n_2}$, splits the $n_1+n_2$ observations into $\sum_{i=0}^2 k_i$
groups of distinct values with respective frequencies $\{n_{j,1}\}
_{j=1}^{k_1}$, $\{n_{\ell,1}\}_{\ell=1}^{k_2}$ and $\{
q_{r,1}+q_{r,2}\}
_{r=1}^{k_0}$. Henceforth, we shall use the shorter notation
\[
\Lambda_{\bm{n},\bm{k}}(\Ac)=\bigl(\tilde\pi_{n_1,n_2},\bm
{X}^{(n_1)},\bm {Y}^{(n_2)}\bigr)^{-1}(
\pi_{n_1,n_2},\Ac)
\]
with $\Ac$ standing for the collection of pairwise disjoint sets $\{
A_{j,1},A_{\ell,2},A_r\dvt j=1,\ldots,k_1; \ell=1,\ldots,k_2;
r=1,\ldots,k_0\}$. Moreover, for any pair of set function $m_1$ and
$m_2$ on $(\X,\Xcr)$ we set $\bm{m}_i^{\bm{n}^{(i)}}(\Ac_i)=\prod_{j=1}^{k_i}m_i^{n_{j,i}}(A_{j,i})$ and $(\bm{m}_1^{\bm
{q}^{(1)}}\times
\bm{m}_2^{\bm{q}^{(2)}})(\Ac_0)=\prod_{r=1}^{k_0}m_1^{q_{r,1}}(A_{r})
m_2^{q_{r,2}}(A_{r})$. By virtue of \eqref{eqpartialexchange} one has
\renewcommand{\theequation}{\arabic{equation}}
\setcounter{equation}{36}
\begin{equation}
\label{eqjointpiobs} \P \bigl[\Lambda_{\bm{n},\bm{k}}(\Ac) \bigr]=\int
_{P_\X^2}\bm {p}_1^{\bm
{n}^{(1)}}(\Ac_1)
\bm{p}_2^{\bm{n}^{(2)}}(\Ac_2) \bigl(
\bm{p}_1^{\bm{q}^{(1)}}\times\bm {p}_2^{\bm{q}^{(2)}}
\bigr) (\Ac_0) \Phi(\ddr p_1,\ddr p_2).
\end{equation}
Since each $\tilde p_i$ is equal, in distribution, to the normalized
measure $\tilde\mu_i/\tilde\mu_i(\X)$ one can proceed in a similar
fashion as in the proof of Proposition~\ref{pr1} and write
\begin{eqnarray*}
\P \bigl[\Lambda_{\bm{n},\bm{k}}(\Ac) \bigr]&=&\frac{1}{\Gamma
(n_1)\Gamma
(n_2)} \int
_0^\infty\,\ddr u \int_0^\infty
\,\ddr v
\\
&&{}\times\E \bigl[\edr^{-u\tilde\mu_1(\X)-v\tilde\mu_2(\X)} \tilde{\bbm{\mu}}_1^{\bm{n}^{(1)}}(
\Ac_1) \tilde{\bbm{\mu }}_2^{\bm
{n}^{(2)}}(
\Ac_2) \bigl(\tilde{\bbm{\mu}}_1^{\bm{q}^{(1)}}\times
\tilde{\bbm{\mu}}_2^{\bm
{q}^{(2)}}\bigr) (\Ac_0) \bigr].
\end{eqnarray*}
Since CRMs give rise to mutually independent random variables when
evaluated on disjoint sets, which identifies the so-called independence
property of CRMs, the expected value in the integral above is shown to
coincide with
\begin{eqnarray*}
&&\E \bigl[\edr^{-u\tilde\mu_1(\X^*)-v\tilde\mu_2(\X^*)} \bigr] \prod_{i=1}^2
\prod_{j=1}^{k_i}\E \bigl[
\edr^{-u\tilde\mu
_1(A_{j,i})-v\tilde\mu_2(A_{j,i})} \tilde\mu_i^{n_{j,i}}(A_{j,i})
\bigr]
\\
&&\quad{}\times\prod_{r=1}^{k_0}\E \bigl[
\edr^{-u\tilde\mu_1(A_r)-v\tilde
\mu_2(A_r)} \tilde\mu_1^{q_{r,1}}(A_r)
\tilde\mu_2^{q_{r,2}}(A_r) \bigr],
\end{eqnarray*}
where $\X^*=\X\setminus\{(\bigcup_{i=1}^2\bigcup_{j=1}^{k_i}
A_{j,i})\cup(\bigcup_{r=1}^{k_0} A_r)\}$.
In the first product, let us consider $i=1$. A~similar line of
reasoning holds for $i=2$ as well. If we set $h_z(u,v)=z(\psi(u)+\psi
(v))+(1-z)\psi(u+v)$, by virtue of the Fa\`a di Bruno formula the
$j$th factor coincides with
\begin{eqnarray*}
&&\E \bigl[\edr^{-u\tilde\mu_1(A_{j,1})-v\tilde\mu_2(A_{j,1})} \tilde\mu_1^{n_{j,1}}(A_{j,1})
\bigr]\\
 &&\quad= (-1)^{n_{j,1}} \frac{\partial^{n_{j,1}}}{\partial u^{n_{j,1}}} \edr^{-G(A_{j,1}) h_z(u,v)}
\\
&&\quad=\edr^{-G(A_{j,1}) h_z(u,v)} \bigl\{G(A_{j,1}) \bigl(z\tau
_{n_{j,1}}(u)+(1-z)\tau_{n_{j,1}}(v)\bigr)+ R_j(A_{j,1})
\bigr\},
\end{eqnarray*}
where $R_j(A_{j,1})$ is a polynomial in $G(A_{j,1})$ of order greater
than $1$ and $G=cP_0$. Moreover, a multivariate version of the Fa\`a di
Bruno formula, see \cite{constsav}, leads to
\begin{eqnarray*}
&&\E \bigl[\edr^{-u\tilde\mu_1(A_{r})-v\tilde\mu_2(A_{r})} \tilde\mu_1^{q_{r,1}}(A_r)
\tilde\mu_2^{q_{r,2}}(A_r) \bigr]
\\
&&\quad= (-1)^{q_{r,1}+q_{r,2}} \frac{\partial^{q_{r,1}+q_{r,2}}}{\partial u^{q_{r,1}}\partial v^{q_{r,2}}} \edr^{-G(A_r)h_z(u,v)}
\\
&&\quad = \edr^{-G(A_r)h_z(u,v)} \bigl\{G(A_r) (1-z)\tau
_{q_{r,1}+q_{r,2}}(u+v)+R_r^*(A_r)\bigr\}
\end{eqnarray*}
with $R_r^*(A_r)$ denoting a polynomial in $G(A_r)$ of degree greater
than $1$. Combining all these facts together, one obtains
\begin{eqnarray*}
&&\P \bigl[\Lambda_{\bm{n},\bm{k}}(\Ac) \bigr]\\
&&\quad=\frac{\prod_{j,\ell
,r}G(A_{j,1})G(A_{\ell,2})
G(A_r)}{\Gamma(n_1)\Gamma(n_2)}
\\
&&\qquad{}\times\sum_{\bm{i}\in\{0,1\}^{k_1}} \sum_{\bm{l}\in\{0,1\}^{k_2}}
(1-z)^{k_0+|\bm{i}|+|\bm{l}|} z^{k_1+k_2-|\bm{i}|-|\bm{l}|}\\
\\
&&\qquad{}\times \int_0^\infty
\int_0^\infty\prod_{j=1}^{k_1}
\tau_{n_{j,1}}(u+i_j v)
\prod_{l=1}^{k_2}
\tau_{n_{l,2}}(\ell_l u+v) \prod_{r=1}^{k_0}
\tau_{q_{r,1}+q_{r,2}} (u+v) \,\ddr\,\ddr v + R_r^{**}(
\Ac),
\end{eqnarray*}
where $R_r^{**}(\Ac)$ is a polynomial of order greater than
$k=k_1+k_2+k_0$ in the variables $G(A_{j,1})$, with $j=1,\ldots,k_1$,
$G(A_{\ell,2})$, with $\ell=1,\ldots,k_2$, and $G(A_r)$, with
$r=1,\ldots,k_0$. It is apparent that the
probability distribution of $(\bm{X}^{(n_1)},\bm{Y}^{(n_2)})$,
conditional on $\tilde\pi_{n_1,n_2}=\pi_{n_1,n_2}$, is absolutely
continuous with respect to $P_0^k$ and recall that $P_0$ is nonatomic.
In order to determine a density of $(\bm{X}^{(n_1)},\bm{Y}^{(n_2)})$,
conditional on $\tilde\pi_{n_1,n_2}=\pi_{n_1,n_2}$, define $\Ac_\ep
$ as
the collection of sets $\{A_{j,1}^\ep,A_{\ell,2}^\ep,A_r^\ep:
j=1,\ldots,k_1; \ell=1,\ldots,k_2; r=1,\ldots,k_0\}$ with
\[
A_{j,1}^\ep\downarrow\{x_j\},\qquad
A_{\ell,2}^\ep\downarrow\{y_\ell\},\qquad
A_r^\ep\downarrow\{z_r\}
\]
as $\ep\downarrow0$. Hence, a version of the conditional density of
$(\bm{X}^{(n_1)},\bm{Y}^{(n_2)})$, conditional on $\tilde\pi
_{n_1,n_2}=\pi_{n_1,n_2}$, with respect to $P_0^k$ and evaluated at
$(\bm{x},\bm{y},\bm{z})$ is proportional to
\[
\lim_{\ep\downarrow0}\frac{\P [\Lambda_{\bm{n},\bm{k}}(\Ac
^\ep
) ]}{
\prod_{j=1}^{k_1}P_0(A_{j,1}^\ep) \prod_{\ell=1}^{k_2}P_0(A_{\ell
,2}^\ep)
\prod_{r=1}^{k_0}P_0(A_r^\ep)}
\]
and, from previous expansion, it can be easily seen to coincide with 1.
And this proves the statement. 

\subsection{Proof of Proposition \texorpdfstring{\protect\ref{pr4}}{4}}\label{appa3}
The probability distribution $\Lc^*$ defined in \eqref{eqlstar} can be
decomposed as follows
\[
\Lc\bigl(\bbm{\theta}^*,\pi_{n_1,n_2},\bbm{\zeta}\bigr) \Lc\bigl(\bm
X^{(n_1)},\bm Y^{(n_2)} | \bbm{\theta}^*, \pi_{n_1,n_2},\bbm{
\zeta}\bigr).
\]
In a similar fashion to the proof of Proposition \ref{pr2}, we use the notation
\[
\Lambda_{\bm{n},\bm{k}}(\Ac)=\bigl(\tilde\pi_{n_1,n_2},\bbm{\theta}^*
\bigr)^{-1}( \pi_{n_1,n_2},\Ac)
\]
with $\Ac$ standing for the collection of pairwise disjoint sets $\{
A_{j,1},A_{\ell,2},A_r\dvt j=1,\ldots,k_1; \ell=1,\ldots,k_2;
r=1,\ldots,k_0\}$. By virtue of \eqref{eqpartialexchange} and by
definition of $\bbm{\zeta}=(\bbm\zeta^{(1)},\bbm\zeta^{(2)})$, one has
%
%
\begin{eqnarray}
\label{eqjointpilat} 
\P \bigl[\Lambda_{\bm{n},\bm{k}}(\Ac) | \bbm
\zeta \bigr]&=& 
\int_{P_{\X}^3}
\bm{p}_1^{\bm{n}^{(1)} \cdot\bbm\zeta_1^*}(\Ac _1)\bm {p}_2^{\bm{n}^{(2)} \cdot\bbm\zeta_2^*}(
\Ac_1)
\nonumber
\\[-8pt]
\\[-8pt]
\nonumber
&&\hspace*{16pt}{}\times\bm{p}_0^{\bm{n}^{(1)} \cdot(\bm1 - \bbm\zeta_1^*)}(\Ac _1)\bm
{p}_0^{\bm{n}^{(2)} \cdot(\bm1 - \bbm\zeta_2^*)}(\Ac_2) \bm {p}_0^{\bm
{q}^*}(
\Ac_0) \Phi'(\ddr p_0, \ddr p_1,
\ddr p_2), 
\end{eqnarray}
where $\Phi'$ corresponds to the probability distribution of the
random vector
\[
\biggl(\frac{\mu_0}{\mu_0(\X)},\frac{\mu_1}{\mu_1(\X)},\frac
{\mu_2}{\mu
_2(\X)} \biggr)
\]
on $P_\X^3$ and we have used vector notation to denote the inner
products $\bm{n}^{(i)}\cdot\bbm{\zeta}_i^*=\sum_{j=1}^{k_i}n_{j,i}\zeta
_{j,i}^*$ and $\bm{n}^{(i)}\cdot(\bm{1}-\bbm{\zeta}_i^*)=\sum_{j=1}^{k_i}n_{j,i}(1-\zeta_{j,i})$ for $i=1,2$. Moreover, note that
\begin{eqnarray*}
\P\bigl[\bbm{\zeta}=\bigl(\bm{a}^{(1)},\bm{a}^{(2)}\bigr) |
\mu_0,\mu_1,\mu_2\bigr]
= \frac{\mu_{1}(\X)^{|\bm{a}^{(1)}|} \mu_{2}(\X)^{|\bm{a}^{(2)}|/2}
\mu_{0}(\X)^{n_1 + n_2 - |\bm{a}^{(1)}|- |\bm{a}^{(2)}|/2}}{(\mu
_0(\X
)+\mu_1(\X))^{n_1} (\mu_0(\X)+\mu_2(\X))^{n_2}}
\end{eqnarray*}
for any $\bm{a}=(\bm{a}^{(1)},\bm{a}^{(2)})\in\{0,1\}^{n_1}\times\{
0,2\}^{n_2}$.
Thus, by similar arguments to those employed in the proofs of
Propositions \ref{pr1} and \ref{pr2}, we can write
\begin{eqnarray*}
\P \bigl[\Lambda_{\bm{n},\bm{k}}(\Ac) , \bbm{\zeta}=\bm{a} \bigr]&=&
\frac
{1}{\Gamma(n_1)\Gamma(n_2)} \int_0^\infty\,\ddr u \int
_0^\infty\,\ddr v u^{n_1-1} v^{n_2-1}
\\
&&\times{}\E \bigl[\edr^{-u(\mu_0(\X)+\mu_1(\X))-v(\mu_0(\X)+\mu_2(\X))} \bbm{\mu}_1^{\bm{n}^{(1)}\cdot\bbm{\zeta}_1^*}(
\Ac_1){\bbm{\mu }}_0^{\bm
{n}^{(1)}\cdot(\bm1 - \bbm{\zeta}_1^*)}(\Ac_1)
\\
&&\hspace*{17pt}{}\times{\bbm{\mu}}_2^{\bm{n}^{(2)}\cdot\bbm\zeta_2^*}(\Ac _2) {\bbm
{\mu}}_0^{\bm{n}^{(2)}\cdot(\bm1 - \bbm\zeta_2^*)}(\Ac_2) \bigl(\bbm{
\mu}_0^{\bm{q}^*}\bigr) (\Ac_0) \bigr],
\end{eqnarray*}
where $\bm{a}$ is a vector such that $\bm{a}^{(i)}$ contains $k_i$
labels $\zeta_{j,i}^*$ such that
\[
\P\bigl[\theta_{j,i}^*\in A | p_0,p_1,p_2
\bigr]=p_{\zeta_{j,i}^*}(A).
\]
Using the independence property of CRMs and the independence of $\mu
_0$, $\mu_1$ and $\mu_2$, the expected value in the integral above can
be rewritten as
\begin{eqnarray*}
&&\E \bigl[\edr^{-u(\mu_0(\X^*)+\mu_1(\X^*))-v(\mu_0(\X^*)+\mu
_2(\X
^*))} \bigr]
\\
&&\quad\times\prod_{i=1}^2 \prod
_{j=1}^{k_i} \E \bigl[\edr^{-u(\mu
_0(A_{j,i})+\mu_1(A_{j,i}))-v(\mu_0(A_{j,i})+\mu_2(A_{j,i}))}\mu
_i(A_{j,i})^{n_{j,i}\zeta_{j,i}^*} \mu_0(A_{j,i})^{n_{j,i}(1-\zeta
_{j,i}^*)}
\bigr]
\\
&&\quad\times\prod_{r=0}^{k_0}\E \bigl[
\edr^{-u(\mu_0(A_{r})+\mu
_1(A_{r}))-v(\mu_0(A_{r})+\mu_2(A_{r}))}\mu _0(A_{r})^{q_{r,1}+q_{r,2}} \bigr],
\end{eqnarray*}
where $\X^*=\X\setminus\{(\bigcup_{i=1}^2\bigcup_{j=1}^{k_i}
A_{j,i})\cup(\bigcup_{r=1}^{k_0} A_r)\}$.
In the first product consider $i=1$, a similar line of reasoning holds
then for $i=2$. The $j$th factor coincides with
%
%
\begin{eqnarray}
\label{jthfactor} &&\E \bigl[\edr^{-v\mu_2(A_{j,1})} \bigr] \E \bigl[
\edr^{-u\mu
_1(A_{j,1})}\mu_1(A_{j,1})^{n_{j,1}\zeta_{j,1}^*} \bigr]
\nonumber
\\[-8pt]
\\[-8pt]
\nonumber
&&\quad{}\times\E \bigl[\edr^{-(u+v)\mu_0(A_{j,1})}\mu _0(A_{j,1})^{n_{j,1}(1-\zeta_{j,1}^*)}
\bigr],
\end{eqnarray}
where
\[
\E \bigl[\edr^{-v\mu_2(A_{j,1})} \bigr]=\edr^{-c P_0(A_{j,1}) \psi(v)}
\]
and, by virtue of the Fa\`a di Bruno formula,
\begin{eqnarray*}
\E \bigl[\edr^{-u\mu_1(A_{j,1})}\mu_1(A_{j,1})^{n_{j,1}\zeta
_{j,1}^*}
\bigr]&=&(-1)^{n_{j,1}\zeta_{j,1}^*}\frac{\partial
^{n_{j,1}\zeta
_{j,1}^*}}{\partial u^{n_{j,1}\zeta_{j,1}^*}} \edr^{-cz P_0(A_{j,1})
\psi(u)}
\\&
=&\edr^{-cz P_0(A_{j,1}) \psi(u)} \bigl\{cz \bigl[P_0(A_{j,1})\tau
_{n_{j,1}}(u)+R_{j,1}(A_{j,1}) \bigr] \bigr
\}^{\zeta_{j,1}^*}
\end{eqnarray*}
and
\begin{eqnarray*}
&&\E \bigl[\edr^{-(u+v)\mu_0(A_{j,1})}\mu_0(A_{j,1})^{n_{j,1}(1-\zeta
_{j,1}^*)}
\bigr]
\\
&&\quad=(-1)^{n_{j,1}(1-\zeta_{j,1}^*)}\frac{\partial
^{n_{j,1}(1-\zeta
_{j,1}^*)}}{\partial s^{n_{j,1}(1-\zeta_{j,1}^*)}} \edr^{-c(1-z)
P_0(A_{j,1}) \psi(s)}\Big\vert
_{s=u+v}
\\
&&\quad=\edr^{-c(1-z) P_0(A_{j,1}) \psi(u+v)} \bigl\{c(1-z) \bigl[P_0(A_{j,1})
\tau_{n_{j,1}}(u+v)+R_{j,1}(A_{j,1}) \bigr] \bigr\}
^{1-\zeta_{j,1}^*}.
\end{eqnarray*}
In the previous expressions, we have agreed that $\partial^0/\partial
s^0$ is the identity operator and that $R_{j,1}(A_{j,1})$ is some
polynomial in $P_0(A_{j,1})$ of order greater than $1$. Thus, the
product in \eqref{jthfactor} is equal to
%
%
\begin{eqnarray}
\label{jthfactor2} \edr^{-c P_0(A_{j,1})h_z(u,v)}c z^{\zeta_{j,1}^*}(1-z)^{1-\zeta
_{j,1}^*}
 \bigl\{P_0(A_{j,1})\tau_{n_{j,1}}\bigl(u+
\bigl(1-\zeta _{j,1}^*\bigr)v\bigr)+R_{j,1}(A_{j,1})
\bigr\}.
\end{eqnarray}
Analogously, one has
%
%
\begin{equation}
\label{firstfactor} \E \bigl[\edr^{-u(\mu_0(\X^*)+\mu_1(\X^*))-v(\mu_0(\X^*)+\mu
_2(\X
^*))} \bigr]=\edr^{-c P_0(\X^*) h_z(u,v)}
\end{equation}
and
%
%
\begin{eqnarray}
\label{rthfactor}&& \E \bigl[\edr^{-u(\mu_0(A_{r})+\mu_1(A_{r}))-v(\mu_0(A_{r})+\mu
_2(A_{r}))}\mu_0(A_{r})^{q_{r,1}+q_{r,2}}
\bigr]
\nonumber
\\[-8pt]
\\[-8pt]
\nonumber
&&\quad=\edr^{-c P_0(A_{r}) h_z(u,v)} c (1-z) \bigl\{P_0(A_r)\tau
_{q_{r,1}+q_{r,2}}(u+v)+R_r(A_r) \bigr\},
\end{eqnarray}
where $R_r(A_r)$ is some polynomial in $P_0(A_r)$ of order greater than $1$.
By combining the expressions \eqref{jthfactor2}--\eqref{rthfactor}, we
obtain that $\P [\Lambda_{\bm{n},\bm{k}}(\Ac) , \bbm{\zeta
}=\bm
{a} ]$ coincides with
\begin{eqnarray*}
&&\frac{c^k z^{|\bbm{\zeta}_{1}^*|+|\bbm{\zeta}_{2}^*|}
(1-z)^{k_1+k_2-|\bbm{\zeta}_{1}^*|-|\bbm{\zeta}_{2}^*|}}{\Gamma
(n_1)\Gamma
(n_2)}P_0^{k}(\Ac)
\\
&&\quad{}\times\int_0^\infty\int_0^\infty
u^{n_1-1} v^{n_2-1}\edr^{-c
h_z(u,v)}\prod
_{j=1}^{k_1}\tau_{n_{j,1}}\bigl(u+\bigl(1-
\zeta_{j,1}^*\bigr)v\bigr)
\\
&&\hspace*{38pt}\qquad{}\times\prod_{j=1}^{k_2}\tau_{n_{j,2}}\bigl(
\bigl(1-\zeta_{j,2}^*\bigr)u+v\bigr) \prod
_{r=1}^{k_0}\tau_{q_{r,1}+q_{r,2}}(u+v)\,\ddr u \,\ddr
v+R^*(\Ac),
\end{eqnarray*}
where $R^{*}(\Ac)$ is a polynomial in the variables $P_0(A_{j,1})$,
with $j=1,\ldots,k_1$, $P_0(A_{\ell,2})$, with $\ell=1,\ldots,k_2$, and
$P_0(A_r)$, with $r=1,\ldots,k_0$, of order greater than
$k=k_1+k_2+k_0$ and $P_0^k(\Ac)=\prod_{i=1}^{k_1}\prod_{j=1}^{k_2}\prod_{r=1}^{k_0}P_0(A_{i,1})P_0(A_{j,2})P_0(A_{r})$.
It is apparent that the probability distribution of $(\bbm\theta
^{(1)},\bbm\theta^{(2)})$, conditional on $\wt{\pi}_{n_1,n_2}=\pi
_{n_1,n_2}$, is degenerate on $\Theta^k$ and the probability
distribution of the distinct values $\bbm{\theta}^*=(\bbm{\theta
}_*^{(1)},\bbm{\theta}_*^{(2)},\bbm{\theta}_*)$ is absolutely continuous
with respect to~$P_0^k$. In order to determine a density of $(\bbm\theta
^*,\bbm\zeta^*,\wt{\pi}_{n_1,n_2})$, introduce $\Ac_\varepsilon$
as in
the proof of Proposition~\ref{pr2} with
\[
A_{j,1}^\ep\downarrow\bigl\{\theta_{j,1}^*\bigr
\},\qquad A_{\ell,2}^\ep\downarrow\bigl\{\theta_{\ell,2}^*\bigr
\},\qquad A_r^\ep\downarrow\bigl\{\theta_r^*\bigr
\}
\]
as $\ep\downarrow0$ and observe that
\[
\lim_{\varepsilon\downarrow0}\frac{\P [\Lambda_{\bm n,\bm
k}(\Ac
_\varepsilon) ]}{P_0^k(\Ac_\varepsilon)}= g\bigl(\bm n^{(1)},\bm
n^{(2)},\bm q^{(1)},\bm q^{(2)},\bbm{\zeta}^*\bigr)
\]
and that
%
%
\begin{equation}
\Lc\bigl(\bbm\theta^*,\bbm\pi_{n_1,n_2},\bbm\zeta\bigr)=g\bigl(\bm
n^{(1)},\bm n^{(2)},\bm q^{(1)},\bm q^{(2)},
\bbm{\zeta}^*\bigr) \prod_{i=1}^k
g_0\bigl(\theta _i^*\bigr). 
\end{equation}
Since the vector $(\bm{X}^{(n_1)},\bm{Y}^{(n_2)})$, given the partition
$\tilde\pi_{n_1.n_2} = \pi_{n_1,n_2}$ and the distinct values $(\bbm
\theta_*^{(1)},\bbm\theta_*^{(2)},\bbm\theta_*)$, is independent
from the
labels $\bbm{\zeta}$, the result follows from \eqref{eqbigprod}.

\subsection{Proof of Corollary \texorpdfstring{\protect\ref{prpjointlab}}{2}}\label{appa4}
If $(\tilde\mu_1,\tilde\mu_2)$ are GM-dependent gamma CRMs, then
one has
$
\tau_q=\Gamma(q)(1+u)^{-q}
$
and
$
\psi(u)=\log(1+u).
$
By plugging these expressions into \eqref{eqconditeppf} and resorting
to identity 3.197.1 in \cite{gr}, we obtain that $g(\bm{n}^{(1)},\bm
{n}^{(2)},\bm{q}^{(1)},\bm{q}^{(2)},\bbm{\zeta}^*)$ is equal to
%
%
\begin{eqnarray}
\label{step1}&& c^k z^{\tilde k_1+\tilde k_2}(1-z)^{\tilde k_0}
\frac{\Gamma
(c+n_1-\bar
n_1)}{\Gamma(n_1)\Gamma(c+n_1+n_2-\bar n_1)} \xi\bigl(\bm{n}^{(1)},\bm {n}^{(2)},\bm{q}^*
\bigr)
\nonumber
\\[-8pt]
\\[-8pt]
\nonumber
&&\quad{}\times\int_0^\infty u^{n_1-1}(1+u)^{-c-n_1+\bar{n}_2}
{}_2F_1(\bar {n}_2+cz,n_2;n_1+n_2-
\bar n_1+c;-u)\,\ddr u,
\end{eqnarray}
where we recall that $\tilde k_0=k_1+k_2-\tilde k_1-\tilde k_2$. The
simple change of variable $t=u/(1+u)$ and the transformation formula
for hypergeometric functions
\[
{}_2F_1(\alpha,\beta;\gamma;z)=(1-z)^{-\alpha}
{}_2F_1\bigl(\alpha,\gamma -\beta ;\gamma;z/(z-1)\bigr)
\]
let us rewrite the integral in \eqref{step1} as
\[
\int_0^1 t^{n_1-1}(1-t)^{c+cz-1}
{}_2F_1(\bar n_2+cz,c+n_1-\bar
n_1; c+n_1+n_2-\bar n_1;t)\,\ddr
t.
\]
The proof is then completed by resorting to identity 7.512.5 in \cite
{gr}. 

\subsection{Proof of Corollary \texorpdfstring{\protect\ref{prpjointlabstable}}{3}}\label{appa5}
If $(\tilde\mu_1,\tilde\mu_2)$ are GM-dependent $\sigma$-stable CRMs,
then one has
$
\tau_q=\sigma(1-\sigma)_{q-1}u^{\sigma-q}
$
and
$
\psi(u)=u^{\sigma}.
$
By plugging these expressions into \eqref{eqconditeppf} we obtain
that $g(\bm{n}^{(1)},\bm{n}^{(2)},\bm{q}^{(1)},\bm{q}^{(2)},\bbm
{\zeta
}^*)$ is equal to
\begin{eqnarray*}
&&\frac{c^k z^{\tilde k_1+\tilde k_2}(1-z)^{\tilde k_0} \sigma
^{k}}{\Gamma(n_1)\Gamma(n_2)} \xi_\sigma\bigl(\bm{n}^{(n_1)},\bm
{n}^{(n_2)},\bm{q}^*\bigr)
\\
&&\quad{}\times\int_0^\infty\int_0^\infty
\frac{u^{n_1-\bar n_1+\tilde k_1
\sigma-1}v^{n_2-\bar n_2+\tilde k_2\sigma-1}(u+v)^{\tilde k_0\sigma
-n_1-n_2+\bar n_1+\bar n_2}}{\exp\{c[z(u^\sigma+v^\sigma
)+(1-z)(u+v)^\sigma]\}}\,\ddr u \,\ddr v.
\end{eqnarray*}
The proof is completed by carefully applying the change of variables
$w=u/(u+v)$ and $s=u+v$.
\end{appendix}

\section*{Acknowledgements} The authors are grateful to an Associate
Editor and three referees for their constructive comments and valuable
suggestions. This work was supported by the European Research Council
(ERC) through StG ``N-BNP'' 306406. Part of the material presented here
is contained in the Ph.D. thesis~\cite{nipoti} defended at the University
of Pavia (Italy) in June 2011.

%



\printhistory

\end{document}